\documentclass[12pt,onecolumn]{IEEEtran}
\usepackage{graphicx}
\usepackage{mathtools,amsmath,amsfonts}
\usepackage{multirow,amssymb}

\usepackage{setspace}%
\usepackage{enumitem,verbatim}
\usepackage{algorithm}
\usepackage{algorithmic}
\usepackage{dsfont}
\usepackage{setspace}
\usepackage[hidelinks]{hyperref}

\DeclareMathOperator*{\argmin}{\arg\!\min}
\DeclareMathOperator*{\argmax}{\arg\!\max}
\usepackage{xcolor}
\usepackage[para,online,flushleft]{threeparttable}
\newtheorem{assumption}{\bf Assumption}
\def\grad{\nabla}
\def\ba{\mathbf{a}}
\def\bb{\mathbf{b}}

\def\bp{\mathbf{p}}
\def\bq{\mathbf{q}}
\def\br{\mathbf{r}}

\def\bD{\mathbf{D}}

\def\bI{\mathbf{I}}

\def\cD{\mathcal{D}}

\def\cF{\mathcal{F}}

\def\cH{\mathcal{H}}

\def\cK{\mathcal{K}}
\def\cL{\mathcal{L}}

\def\cN{\mathcal{N}}
\def\cO{\mathcal{O}}
\def\cP{\mathcal{P}}

\def\cU{\mathcal{U}}

\def\cX{\mathcal{X}}
\def\cY{\mathcal{Y}}
\def\cZ{\mathcal{Z}}

\def\mE{\mathbb{E}}

\def\abs#1{\left|#1\right|}
\def\norm#1{\left\|#1\right\|}

\newcommand{\reals}{\mathbb{R}}
\newcommand{\integers}{\mathbb{Z}}

\newcommand{\dom}{\mathop{\bf dom}}

\newcommand{\diag}{\mathop{\bf diag}}

\def\fprod#1{\left\langle#1\right\rangle}

\def\ind#1{\mathbb{I}_{#1}}

\def\ones{\mathbf{1}}
\def\red#1{\textcolor{black}{#1}}
\newcommand{\intr}{\mathop{\bf int}}
\def\eyh#1{\textcolor{black}{#1}}

\newtheorem{theorem}{\bf Theorem}[section]
\newtheorem{corollary}{\bf Corollary}[section]
\newtheorem{lemma}[theorem]{\bf Lemma}
\newtheorem{definition}{\bf Definition}[section]
\newtheorem{remark}{\bf Remark}[section]

\newcommand{\qed}{\mbox{}\hspace*{\fill}\nolinebreak\mbox{$\rule{0.7em}{0.7em}$}}

\begin{document}

\title{A Stochastic Variance-reduced Accelerated Primal-dual Method for Finite-sum Saddle-point Problems}

%\subtitle{Variance-reduced Accelerated Primal-dual Method}
%\titlerunning{A Stochastic Variance-reduced Accelerated Primal-dual Method}  
%\author{Erfan Yazdandoost Hamedani \and Afrooz Jalilzadeh }

\author{Erfan Yazdandoost Hamedani, Afrooz Jalilzadeh \\ {\small
             Department of Systems and Industrial Engineering, The University of Arizona
             
              erfany@arizona.edu, afrooz@arizona.edu }
           %\and
            %Afrooz Jalilzadeh,  Corresponding author  \at
             %  Department of Systems and Industrial Engineering, The University of Arizona, 1127 E James Rogers Way, Tucson, AZ\\
              %afrooz@arizona.edu
}

%The correct dates will be entered by the editor.

\maketitle
\allowdisplaybreaks
\begin{abstract}
In this paper, we propose a variance-reduced primal-dual algorithm \eyh{with Bregman distance} for solving convex-concave saddle-point problems with finite-sum structure and nonbilinear coupling function. 
This type of problems typically arises in machine learning and game theory.
Based on some standard assumptions, the algorithm is proved to converge with oracle complexity of $\cO(\frac{\sqrt{n}}{\epsilon})$  and $\cO(\frac{n}{\sqrt{\epsilon}}+\frac{1}{\epsilon^{1.5}})$ using constant and non-constant parameters, respectively where $n$ is the number of function components.
Compared with existing methods, our framework yields a significant improvement over the number of required primal-dual gradient samples to achieve $\epsilon$-accuracy of the primal-dual gap. We tested our method for solving a distributionally robust optimization problem to show the effectiveness of the algorithm. 
\end{abstract}
%\keywords{Primal-dual\and Saddle-point\and first-order method\and Variance reduction}

%All acknowledgements should be placed in the back of the paper after Conclusions..

\section{Introduction}
\label{sec:intro}
Let $(\cX,\norm{\cdot}_\cX)$ and $(\cY,\norm{\cdot}_\cY)$ be finite-dimensional normed vector spaces, with dual spaces $(\cX^*,\norm{\cdot}_{\cX^*})$ and $(\cY^*,\norm{\cdot}_{\cY^*})$, respectively, and $\cN\triangleq \{1,\hdots,n\}$. We consider the following convex-concave saddle-point (SP) problem
\begin{align}\label{eq:sp-problem}
\min_{x\in X}\max_{y\in Y} \cL(x,y)\triangleq f(x)+\Phi(x,y)-h(y),
\end{align}
where $X\subseteq \cX$ and $Y\subseteq \cY$ are nonempty, closed, and convex sets; $f:\cX\to \bar{\reals}\triangleq \reals\bigcup\{+\infty\}$ and $h:\cY\to\bar{\reals}$ are convex, closed, and proper functions (possibly nonsmooth); moreover, $\Phi(x,y)\triangleq \frac{1}{n}\sum_{i\in\cN}\Phi_i(x,y)$ is  a convex-concave function, i.e., $\Phi(\cdot,y)$ is  convex for any $y\in\cY$ and $\Phi(x,\cdot)$ is concave for any $x\in\cX$, and $\Phi_i$ satisfies certain differentiablity assumptions for $i\in\cN$ -- see Assumption \ref{assum:lip}. 

We are motivated by designing an efficient algorithm to solve \eqref{eq:sp-problem} which emerges in machine learning and data analysis problems when $n$ is large. There has been a lot of efforts to solve large-scale optimization problems efficiently using different approaches, e.g., variance reduction and block-coordinate schemes. Specifically, when the objective has a finite-sum structure, variance reduction schemes close the oracle complexity gap between the deterministic and stochastic settings by  providing an unbiased estimator of gradients reducing the variance of the error of gradient estimator, e.g., SAG \cite{roux2012stochastic}, SVRG \cite{johnson2013accelerating}, SAGA \cite{defazio2014saga}. On the other hand, with the emerging complexities arising in different areas, SP problems are becoming more popular, and different methods have been introduced to solve such problems. Unlike optimization problems, the number of methods solving large-scale SP problems with finite-sum structure is limited, most of which only consider the strongly-convex strongly-concave setting. Our goal in this paper is to introduce an SVRG-type variance reduction technique for primal-dual algorithms \eyh{with Bregman distance} for solving convex-concave SP problem \eqref{eq:sp-problem}. 

{\bf Notations.}
Let $\mathbb{S}^n_{++}$ ($\mathbb{S}^n_+$) be the set of $n\times n$ symmetric positive (semi-) definite matrices, $\bI_n$ {denotes} the $n\times n$ identity matrix, and $\ones_n$ denotes an $n$-dimensional vector of ones. $\mE[\cdot]$ denotes the expectation operation and $\widetilde{\cO}(\cdot)$ denotes $\cO(\cdot)$ up to a logarithmic factor. %In the rest, $f(\cdot)=\Theta(g(\cdot))$ means $f(\cdot)=\Omega(g(\cdot))$ and $f(\cdot)=\cO(g(\cdot))$.
\vspace{-5mm}
\subsection{Applications} \label{sec:application}
There is a wide range of real-life problems arising in machine learning, image processing, game theory, etc. such that they can be formulated as a special case of \eqref{eq:sp-problem}. We briefly introduce some of the interesting examples below. 

{\bf I. Distributionally robust optimization (DRO):} Let $(\Omega,\cF,\mathbb{P})$ be a probability space where $\Omega=\{\zeta_1,\hdots, \zeta_n\}$, $\ell:X\times \Omega\to \reals$ is a convex loss function, and we define $\ell_i(u)\triangleq \ell(u;\zeta_i)$. DRO studies worse case performance under uncertainty to find solutions with some specific confidence level \cite{namkoong2016stochastic}. This problem can be formulated as follows:
\begin{align}\label{eq:DRO}
    \min_{u\in X}\max_{y\in\cP}~\mE_{\zeta\sim \mathbb{P}}[\ell(u;\zeta)]=\sum_{i=1}^n y_i\ell_i(u) ,
\end{align}
where $\cP$ represents the uncertainty set.
For instance, $\cP=\{y\in\Delta_n : V(y,\frac{1}{n}\ones_n)\leq \rho\}$ is an uncertainty set considered in different papers such as \cite{namkoong2016stochastic}, where $\Delta_n\triangleq \{y=[y_i]_{i=1}^n\in\reals^n_+ \mid \sum_{i=1}^n y_i=1\}$ is an $n$-dimensional simplex set, and $V(Q,P)$ denotes the divergence measure between two sets of probability measures $Q$ and $P$. Using a variable $\lambda\in\reals_+$ we can relax the divergence constraint in \eqref{eq:DRO} to obtain the following equivalent problem:
\begin{align}\label{eq:DRO2}
    \min_{\substack{u\in X\\ \lambda\geq 0}}\max_{\substack{y=[y_i]_{i=1}^n, \\ y\in\Delta_n}}~\sum_{i=1}^n y_i\ell_i(u) - \frac{\lambda}{n} \big(V(y,\tfrac{1}{n}\ones_n)- \rho\big) ,
\end{align}
Let $x=[u^\top;\lambda]^\top$, \eqref{eq:DRO2} is a special case of \eqref{eq:sp-problem} by defining $\Phi_i(x,y)\triangleq ny_i\ell_i(u)-\lambda(\frac{1}{2}(ny_i-1)^2-\frac{\rho}{n})$, $f(x)=\ind{X\times \reals_+}(u,\lambda)$, and $h(y)=\ind{\Delta_n}(y)$.

{\bf II. Learning a kernel matrix:} Suppose we are given a set of labeled data points consisting of feature vectors $\{\ba_i\}_{i=1}^n\subset\reals^m$, and the corresponding labels $\{b_i\}_{i=1}^n\subset\{-1,+1\}$. Consider $N$ different embedding of the data and let  $K_i\in\mathbb{S}^n_{+}$ be the corresponding kernel matrix. The objective is to learn a kernel matrix $K$ belonging to a class of kernel matrices which is a convex set generated by $\{K_i\}_{i=1}^N$, i.e. $K\in\cK\triangleq \{\sum_{i=1}^M y_i K_i:\ y_i\geq 0,\ i=1,\ldots, M\}$, such that it minimizes the training error of a kernel SVM as a function of $K$ -- see~\cite{lanckriet2004learning} for more details. Then one needs to solve the following problem:
\begin{align}
\label{eq:kernel_learn_simple}
\min_{\substack{y\in\reals^M_+\\ \fprod{\br,y}=c,}}\ \max_{\substack{x:\ 0\leq x\leq C\ones_n, \\ \ \ \fprod{\bb,x}=0}} 2x^\top\ones_n - \sum_{i=1}^M {y_i}x^\top H(K_i)x-\lambda\norm{x}_2^2,
\end{align}
where $c,C\geq 0$ and $\lambda\geq 0$  are model parameters, $y=[y_i]_{i=1}^n$, $\br=[\text{trace}(K_i)]_{i=1}^M$, $\bb=[b_i]_{i=1}^M$ and $H(K_i)\triangleq \diag(\bb)K_i\diag(\bb)$. Clearly \eqref{eq:kernel_learn_simple} has a finite sum objective and is a special case of \eqref{eq:sp-problem}.

{\bf III. Two-player zero-sum game with a nonlinear payoff:} This problem arising in game theory, considers computing the equilibrium of a convex-concave two-player game of the following form:
\begin{align}\label{game}
    \min_{\substack{\bq=[q_i]_{i=1}^N\\\fprod{\ones_N,\bq}=Q}}\max_{\substack{\bp=[p_i]_{i=1}^N\\ \fprod{\ones_N,\bp}=P}} &\sum_{i=1}^N \log\Big(1+\frac{\beta_ip_i}{\sigma_i+q_i}\Big)~+\sum_{j=1}^M \log(1+\exp(\fprod{u_j,\bq})) \\
    &-\sum_{j=1}^M \log(1+\exp(\fprod{v_j,\bp})). \nonumber
\end{align}
The problem in \eqref{game} includes some interesting special cases such as the water filling  problem arising in information theory (see \cite{boyd2004convex,chen2017accelerated}) when $u_j=v_j=0$, for all $j\in\{1,\hdots,M\}$. In particular, consider $N$ given Gaussian communication channels each having signal power $p_i$ and noise power $q_i$, for $i\in\{1,\hdots,N\}$. From Shannon-Hartley equation, the maximum capacity of channel $i$ is proportional to $\log(1 +\frac{\beta_ip_i}{\sigma_i+q_i} )$ where $\sigma_i$ is the receiver noise and $\beta_i>0$ is a constant. The goal is to maximize
the total capacity given total power $P > 0$ while an adversary aims to reduce the total capacity given total noise power $Q > 0$. Therefore, we aim to allocate the signal power such that with the worst allocation of noise power, we obtain
the largest total capacity for the system. 
% which can be formulated as follows:
% \begin{align}
%     \min_{\substack{q=[q_i]_{i=1}^N\\\fprod{\ones,q}=Q}}\max_{\substack{p=[p_i]_{i=1}^N\\ \fprod{\ones,p}=P}} \sum_{i=1}^N \log\Big(1+\frac{\beta_ip_i}{\sigma_i+q_i}\Big)
% \end{align}

\subsection{Related work}  
\eyh{SP problems have become increasingly popular in recent years due to their ability to solve a wider range of problems.} There have been several \eyh{studies} on deterministic first-order primal-dual algorithms for solving \eqref{eq:sp-problem} when $\Phi(x,y)$ is bilinear, i.e., $\Phi(x,y)=\fprod{Ax,y}$, such as \cite{dang2014randomized,chambolle2016ergodic,he2016accelerated,wang2017exploiting,du2018linear}, and few others have considered a more general non-bilinear setting \cite{palaniappan2016stochastic,nemirovski2004prox,juditsky2011first,he2015mirror,kolossoski2017accelerated} \eyh{in which an} optimal rate of $\cO(\frac{1}{\epsilon})$ has been shown for the convex-concave setting. However, when considering problem \eqref{eq:sp-problem}, this rate is directly affected by the number of function components; hence, the oracle complexity (number of primal-dual sample gradients) is $\cO(\frac{n}{\epsilon})$ which requires a high computational effort for large-scale problems, i.e., $n$ is large.

Many stochastic primal-dual algorithms have been introduced in different studies with the aim of addressing more general problems having an expectation in the objective function, and achieving lower per iteration complexity. However, this comes at the cost of dropping the oracle complexity to $\cO(\frac{1}{\epsilon^2})$ for convex-concave setting -- see \cite{nemirovski2009robust,juditsky2011solving,zhao2019optimal}. After introducing variance reduction techniques for optimization problems, different attempts have been made to adopt such techniques in primal-dual algorithms for solving \eqref{eq:sp-problem} or its special case when $\Phi$ is linear in $y$; however, most of these studies only focused on strongly-convex strongly-concave setting such as \cite{palaniappan2016stochastic,lian2017finite,zhang2019composite,devraj2019stochastic}, and few others \cite{hien2017inexact,yan2019stochastic} consider a more general setting which we will briefly describe next.  

In \cite{hien2017inexact}, a randomized primal-dual smoothing technique has been proposed \eyh{for solving \eqref{eq:sp-problem} by inexactly solving a sequence of subproblems when the $\Phi(x,y)$ has a finite-sum structure. 
Assuming that each $\Phi_i(x,y)$ is smooth and convex-concave, and $f$ is strongly convex primal and dual oracle complexities of $\widetilde{\cO}((n+\sqrt{n\kappa})/\sqrt{\epsilon})$ and $\widetilde{\cO}(\sqrt{n}/\epsilon)$, where $\kappa$ denotes the condition number, have been shown, respectively.}
%the objective is strongly convex in $x$, and $X$ is bounded, an oracle complexity of $\widetilde{\cO}(\sqrt{n}/\epsilon)$ has been shown. However, in this work we show that our method can achieve ${\cO}(\sqrt{n}/\epsilon)$ without assuming strong convexity and boundedness of $X$. 
In \cite{yan2019stochastic} the problem of $\min_{x\in X}\max_{y\in Y} f(x)+y^\top g(x)-h(x)$ has been considered. %although they need not to assume the proximal function of $f$ and $h$ can be computed efficiently; hence, they use subgradients of nonsmooth functions $f,h$
This problem can be also equivalently written as $\min_{x\in X} P(x)$ where $P(x) =h^*(g(x))+f(x)$ and $h^*(\cdot)$ denotes the convex conjugate of function $h(\cdot)$. They proposed a restarted stochastic primal-dual algorithm (RSPD) in which noisy partial primal and dual (sub)gradients are used. The algorithm is restarted periodically and in the outer loop $\argmax_{y\in\dom h^*}\cL(x,y)$ for some given $x$ is required to be computed efficiently; however, this operation might be computationally expensive for a more general problem \eqref{eq:sp-problem}. Assuming that the partial (sub)gradients are bounded, $h^*$ follows the H\"older condition with constants $(L,\nu)$, $g(\cdot)$ is Lipschitz continuous, and $P(\cdot)$ satisfies lower error bound with parameter $\theta\geq 0$, i.e., $\hbox{dist}(x,X^*)\leq c(P(x)-P(x^*))^\theta$ for some $c>0$ where $X^*$ denotes the optimal set, they demonstrated an oracle complexity of $\cO(\frac{1}{\epsilon^{2(1-\nu\theta)}})$. Note that if $P(\cdot)$ does not obey the lower error bound condition, i.e., $\theta=0$, then their oracle complexity is $\cO(\frac{1}{\epsilon^2})$. 
In our recent study \cite{jalilzadeh2019randomized} we consider problem \eqref{eq:sp-problem} where $x$ and $y$ are assumed to have $M$ and $N$ blocks, respectively and the problem has a coordinate-friendly structure. A doubly stochastic block-coordinate primal-dual algorithm has been proposed with single and increasing batch-size in which at each iteration only one block of $x$ and $y$ are updated. Assuming that the partial gradients of $\Phi$ are bounded, oracle complexities of $\widetilde{\cO}(\frac{1}{\epsilon^2})$ for single and $\widetilde{\cO}(\frac{1}{\epsilon}\min\{n,\frac{1}{\epsilon}\})$ for increasing batch-size have been achieved.  
%\footnote{$\widetilde{\cO}(\cdot)$ denotes ${\cO}(\cdot)$ up to a logarithmic factor.} 

In contrast to the existing studies mentioned above, we aim to \eyh{obtain an improved oracle complexity under weaker assumptions} for problem \eqref{eq:sp-problem} where $\Phi$ is convex-concave and is neither linear in $x$ nor in $y$.

\subsection{Contribution}
We study SP problem \eqref{eq:sp-problem} with a finite-sum structure where the coupling function $\Phi$ is not linear in $x$ nor in $y$. We develop a stochastic variance-reduced accelerated primal-dual algorithm (SVR-APD) \eyh{with Bregman distance} which is a novel SVRG-type primal-dual algorithm, for solving this problem. Our idea is to consider a new momentum which is a convex combination of the current iterate point and the average of past iterates. This idea combined with an acceleration in terms of partial gradients of $\Phi$ leads to convergence guarantees in terms of the \eyh{standard} gap function $\mE[\sup_{(x,y)\in X\times Y} \cL(x^{(K)},y)-\cL(x,y^{(K)})]$ where $x^{(K)}$ and $y^{(K)}$ denote the ergodic average of the iterates. More precisely, we demonstrate the oracle complexities of $\cO(\frac{\sqrt{n}}{\epsilon})$ and $\cO(\frac{n}{\sqrt{\epsilon}}+\frac{1}{\epsilon^{1.5}})$ using constant and non-constant parameters, respectively. 

Comparing to deterministic methods, our oracle complexity of $\cO(\frac{\sqrt{n}}{\epsilon})$ shows a clear improvement in the order of $\sqrt{n}$ magnitude and it has a lower complexity in comparison with $\cO(1/\epsilon^2)$ for stochastic methods when $\epsilon\leq \mathcal O(\frac{1}{\sqrt{n}})$ -- see Table \ref{tab:my_label}. 
%For instance, if $n=10^8$, then for $\epsilon< 10^{-4}$ our algorithm can achieve an $\epsilon$-accuracy with lower oracle complexity in compare with stochastic methods. 
Moreover, our second result with non-constant parameters leads to an oracle complexity of $\cO\left(\frac{n}{\sqrt{\epsilon}}+\frac{1}{\epsilon^{1.5}}\right)$, which is lower than the stochastic schemes and for $\epsilon\geq \mathcal O({1\over n^2})$ is lower than \eyh{the deterministic counterpart}. \eyh{Finally, comparing} our results, selecting constant parameters leads to a better complexity for a medium to high accuracy of $\epsilon\leq \mathcal O({1\over n})$. 

Furthermore, we were able to incorporate Bregman-distance functions in an SVRG-type method. 
%under some mild conditions -- see Assumption \ref{assum:bregman}.
This can  significantly improve the applicability of the algorithm for some specific problems, e.g., when the constraint set is a simplex set a closed form solution for the projection can be computed using the entropy-distance function rather than projecting onto the set in the Euclidean space. 

%To the best of our knowledge, this is the first result that improves the oracle complexity for general convex-concave SP problems. 

\begin{table}[htb]
    \centering
    \begin{tabular}{|c|c|c|}
    \hline
      Method & Reference &  Oracle complexity\\ \hline
        Deterministic & \cite{nemirovski2004prox}, \cite{he2015mirror}, \cite{malitsky2018proximal},  \cite{hamedani2018primal} & $\cO\left(n/{\epsilon}\right)$\\ \hline
        Stochastic & \cite{nemirovski2009robust}, \cite{juditsky2011solving}, \cite{zhao2019optimal}& $\cO\left(1/{\epsilon^2}\right)$\\ \hline
         \multirow{2}{*}{SVR-APD} & Corollary \ref{cor:constant} & $\cO(\sqrt{n}/\epsilon)$ \\\cline{2-3}   (this paper)& Corollary \ref{cor:nonconstant} &  $\cO\left(\frac{n}{\sqrt{\epsilon}}+\frac{1}{\epsilon^{1.5}}\right)$\\ \hline
    \end{tabular}
    \caption{Comparison of oracle complexity (number of primal-dual gradient samples) between the existing methods for solving SP problem \eqref{eq:sp-problem}.}
    \label{tab:my_label}
\end{table}
\vspace{-5mm}
\subsection{Organization of the Paper}
In the next section, we precisely state our assumptions, describe the proposed SVRG-type  algorithm, and present the oracle complexities of the method under different choices of parameters which are the main results of this paper. Subsequently, in Section~\ref{sec:analysis}, we provide a convergence analysis proving the main results. Later, in Section~\ref{sec:num}, %a numerical example is provided in which
we apply our SVR-APD method to solve the DRO problem and compare it with competitive methods. %Finally, Section~\ref{sec:conclude} concludes the paper.
%\eyh{Moreover, some well-known results are presented in Appendix (section \ref{append}).}
\section{Proposed Method}
For the optimization problem of $\min_{x\in X}f(x)+\frac{1}{n}\sum_{i=1}^n g_i(x)$,  variants of SVRG method have been developed in \cite{allen2016improved} when the objective function is merely convex. The main idea is to keep a full gradient at $\tilde{x}^{k-1}$ in the outer loop and use it to provide an estimate of the full gradient, $\xi_t$, such that a tight upper bound on $\norm{\xi_t-\frac{1}{n}\sum_{i=1}^n\grad g_i(\tilde{x}^{k-1})}$ can be obtained -- see \cite[Lemma A.2]{allen2016improved}. However, such an upper bound cannot be obtained for general saddle-point problems. In this section, we propose the Stochastic Variance-reduced Accelerated Primal-dual (SVR-APD) algorithm displayed in Algorithm \ref{alg:SVRPDA}. %The algorithm comprises two loops; in the outer loop a full primal-dual gradient is computed which is used to update the iterates in the inner loop. The inner loop comprises an accelerated dual step followed by a primal step. 
Our novel idea to resolve this issue is to consider a combination of an iterate with an average of the last loop (see lines 9 and 14) which will help us  deriving some upper bounds for the error of gradient estimates -- see Lemma \ref{lem:upper-noise}.

Next, we provide some fundamental definitions and state our assumptions.

\begin{definition}\label{def:bregman}
Let $\psi_\cX:\cX\to \reals$ and $\psi_\cY:\cY\to \reals$ be continuously differentiable functions on $\intr(\dom f)$ and $\intr(\dom h)$, respectively. Moreover, $\psi_\cX$ and $\psi_\cY$ are 1-strongly convex with respect to $\norm{\cdot}_\cX$ and $\norm{\cdot}_\cY$, respectively. We define the \emph{Bregman distance} function corresponding to the distance generating function $\psi_\cX$ as $\bD_X(x,\bar{x})\triangleq \psi_\cX(x)-\psi_\cX(\bar{x})-\fprod{\grad\psi_\cX(\bar{x}),x-\bar{x}}$, for all $x\in X$ and $\bar{x}\in X^\circ\triangleq X\cap\intr(\dom f)$. Similarly we define $\bD_Y(y,\bar{y})\triangleq \psi_\cY(y)-\psi_\cY(\bar{y})-\fprod{\grad\psi_\cY(\bar{y}),y-\bar{y}}$, for all $y\in Y$ and $\bar{y}\in Y^\circ\triangleq  Y\cap \intr(\dom h)$. Moreover, we define the Bregman diamateres of $X$ and $Y$ under $\psi_\cX$ and $\psi_\cY$ as $B_X$ and $B_Y$, respectively, i.e., $B_X\triangleq \sup_{x\in X,\bar{x}\in X^\circ}\bD_X(x,\bar{x})$ and $B_Y\triangleq \sup_{y\in Y,\bar{y}\in Y^\circ}\bD_Y(y,\bar{y})$.
\end{definition}

\begin{assumption}\label{assum:lip} Let $\bD_{X}$ and $\bD_Y$ be some Bregman distance functions as in Definition~\ref{def:bregman}; $f$ and $h$ are closed convex functions; %functions,
and $\Phi$ is continuously differentiable %function
such that\\
{\bf (i)} for any %fixed
$y\in\dom h\subset \cY$, $\Phi(\cdot,y)$ is convex; for any $i\in\cN$ $\Phi_i(\cdot,y)$ is differentiable; there exist $L_{xx}\geq 0$ and $L_{xy}>0$ such that for all $x,\bar{x}\in\dom f$ and $y,\bar{y}\in\dom h$, one has
%such that for any $x,\bar{x}\in\cX$,
\begin{equation}
\label{eq:Lxx}
\norm{\grad_x \Phi_i(x,y)-\grad_x \Phi_i(\bar{x},\bar{y})}_{\cX^*}\leq L_{xx}\norm{x-\bar{x}}_\cX+L_{xy}\norm{y-\bar{y}}_{\cY},\ \forall i\in\cN,
\end{equation}
{\bf (ii)} for any %fixed
$x\in\dom f\subset \cX$, {$\Phi(x,\cdot)$} is concave; for any $i\in\cN$ $\Phi_i(x,\cdot)$ is differentiable; %in $y$;
{there exist $L_{yx}>0$ and $L_{yy}\geq 0$ such that for all $x,\bar{x}\in\dom f$ and $y,\bar{y}\in\dom h$, one has}
%such that for any $y,\bar{y}\in\cY$ {and $x,\bar{x}\in\cX$,}
%\begin{subequations}
\begin{align}
\label{eq:Lipschitz_y}
\norm{\grad_y \Phi_i(x,y)-\grad_y \Phi_i(\bar{x},\bar{y})}_{\cY^*}\leq L_{yy}\norm{y-\bar{y}}_\cY+{L_{yx}}\norm{x-\bar{x}}_\cX, \ \forall i\in\cN.
%\quad \forall y,\bar{y}\in\cY, \forall x,\bar{x}\in\cX.
%&\norm{\grad_y \Phi(x,y)-\grad_y \Phi(\bar{x},{y})}_{\cY^*}\leq {L_{yx}}\norm{x-\bar{x}}_\cX.
\end{align}
\end{assumption}

Note that \eqref{eq:Lxx} and convexity of $\Phi(\cdot,y)$ imply that for any $y\in\dom h$ and  $x,\bar{x}\in\dom f$,
\begin{align}
0 &\leq \Phi(x,y) - \Phi(\bar{x},y)-\fprod{\grad_x\Phi(\bar{x},y),~x-\bar{x}} \leq\frac{L_{xx}}{2}\norm{x-\bar{x}}_\cX^2.
%,\quad \forall x,\bar{x}\in\dom f.
%\subset \cX.
\label{eq:Lxx_bound}
\end{align}%

\begin{algorithm}
\caption{SVR-APD}
\label{alg:SVRPDA}
\begin{onehalfspacing}
\begin{algorithmic}[1]
\STATE Initialization: $\tilde{x}^0\in\cX,\tilde{y}^0\in\cY$
%\STATE $(\tilde{x}^{0},\tilde{y}^{0})\gets (\tilde{x}^1,\tilde{y}^1$)
\STATE \eyh{$(x_0^1,y_0^1)\gets (\tilde x^0,\tilde y^0)$}
\STATE \eyh{$(r^0,s^0)\gets (\grad \psi_\cX(x_0^1),\grad \psi_\cY(y_0^1))$}
\FOR{$k=1,\hdots,K$}
\STATE $G_y \gets \grad_y\Phi(\tilde x^{k-1},\tilde y^{k-1})$
\STATE $G_x \gets \grad_x\Phi(\tilde x^{k-1},\tilde y^{k-1})$
%\STATE $(x^k_{-1},y^k_{-1})\gets (x^{k-1}_{T^{k-1}-1},y^{k-1}_{T^{k-1}-1})$ and $(x^k_0,y^k_0)\gets (x^{k-1}_{T^{k-1}},y^{k-1}_{T^{k-1}})$
\FOR{$t=0,\hdots,T^k-1$} 
\STATE Pick $j_k\in\cN$ uniformly at random
%\STATE $\hat y_t^k\gets y_t^k-\gamma_y^k(y_t^k-\tilde y^{k-1})$
\STATE $\hat y_t^k\gets \grad\psi_\cY^*\big((1-\gamma_y^k)\grad\psi_\cY(y_t^k)+\gamma_y^k \eyh{s^{k-1}}\big)$
%\STATE {\bf Option I:} $\xi_t\gets \grad_y\Phi_j(x_t^k,y_t^k)-\grad_y\Phi_j(x_t^k,\tilde y^{k-1})+\grad_y\Phi(x_t^k,\tilde y^{k-1})$
\STATE $\xi_t\gets \grad_y\Phi_{j_k}(x_t^k,y_t^k)-\grad_y\Phi_{j_k}(\tilde x^{k-1},\tilde y^{k-1})+G_y$
\STATE $q_t^k\gets \grad_y\Phi_{j_k}(x_t^k,y_t^k)-\grad_y\Phi_{j_k}(x_{t-1}^k,y_{t-1}^k)$
\STATE $y_{t+1}^k\gets \argmin_{y\in Y}\{h(y)-\fprod{\xi_t+  q_t^k,y}+\frac{1}{\sigma^k}\bD_Y(y,\hat{y}_t^k)\}$%(\hat y_t^k+\sigma^k(\xi_t+q_t^k))$ %\red{*$\xi_{t-1}$ has an issue should be explained!}
\STATE Pick $i_k\in\cN$ uniformly at random
%\STATE $\hat x_t^k\gets x_t^k-\gamma_x^k(x_t^k-\tilde x^{k-1})$
\STATE $\hat x_t^k\gets \grad \psi_\cX^*\big((1-\gamma_x^k)\grad\psi_\cX(x_t^k)+\gamma_x^k \eyh{r^{k-1}}\big)$
%\STATE {\bf Option I:} $\zeta_t\gets \grad_x\Phi_i(x_t^k,y_{t+1}^k)-\grad_x\Phi_i(\tilde x^{k-1},y_{t+1}^k)+\grad_x\Phi(\tilde x^{k-1},y_{t+1}^k)$
\STATE $\zeta_t\gets \grad_x\Phi_{i_k}(x_t^k,y_{t+1}^k)-\grad_x\Phi_{i_k}(\tilde x^{k-1},\tilde{y}^{k-1})+G_x$
\STATE $x_{t+1}^k\gets \argmin_{x\in X}\{f(x)+\fprod{\zeta_t,x}+\frac{1}{\tau^k}\bD_X(x,\hat{x}_t^k)\}$%\prox{\tau^k f}(\hat x_t^k-\tau^k \zeta_t)$
\ENDFOR
\STATE $(\tilde x^k,\tilde y^k) \gets \frac{1}{T^k}\sum_{t=0}^{T^k-1} (x_{t+1}^k,~ y_{t+1}^k)$
\STATE \eyh{$(r^k,s^k)\gets \frac{1}{T^k}\sum_{t=0}^{T^k-1} (\grad\psi_\cX(x_{t+1}^k),\grad\psi_\cY(y_{t+1}^k))$}
\STATE \eyh{$(x_0^{k+1},y_0^{k+1},x_{-1}^{k+1},y_{-1}^{k+1})\gets (x_{T^k}^k,y_{T^k}^k,x_{T^k-1}^k,y_{T^k-1}^k)$}
\ENDFOR
\end{algorithmic}
\end{onehalfspacing}
\end{algorithm}

\begin{remark}
It is worth emphasizing that we have only assumed $\Phi(\cdot, \cdot)$ to be a convex-concave function while each component function $\Phi_i(\cdot,\cdot)$ may not be. 
%Moreover, we assumed that the Lipschitz constants of $\grad_x\Phi_i$ and $\grad_y\Phi_i$ are independent of $i$ to facilitate the analysis without loss of generality. Note that, if for any $i\in\{1,\hdots,n\}$ the Lipschitz constants for $\grad_x\Phi_i$ are $L_{xx}^i$ and $L_{xy}^i$, then the Lipschitz constants for $\grad_x\Phi$ will be the average of those constants, i.e., $L_{xx}=\frac{1}{n}\sum_{i=1}^nL_{xx}^i$ and $L_{xy}=\frac{1}{n}\sum_{i=1}^nL_{xy}^i$, respectively. A similar argument can be presented for $L_{yy}^i$ and $L_{yx}^i$.
\end{remark}

\begin{assumption}\label{assum:bregman}
\eyh{The Bregman diameters $B_X$ and $B_Y$ are bounded.}
%Consider the Bregman distance functions defined in Definition \ref{def:bregman}.
%\begin{itemize}
    %\item[a)] The Bregman diameters $B_X$ and $B_Y$ are bounded. 
%    \item[a)] $\psi_\cX$ and $\psi_\cY$ are continuously  differential on $\intr(\dom f)$ and $\intr(\dom h)$, respectively.
%    \item[b)] For any $x\in X$ and $y\in Y$,  $\bD_X(x,\cdot)$ and $\bD_Y(y,\cdot)$ are convex functions.
%    \item[c)] \eyh{The Bregman diameters $B_X$ and $B_Y$ are bounded.}
%\end{itemize}
\end{assumption}

%\begin{remark}
%%For finite dimensional real vector spaces, Assumption \ref{assum:bregman}-$(a)$ is equivalent to the boundedness of $X$ and $Y$ in many scenarios such as when $\psi_\cX(\cdot)=\frac{1}{2}\norm{\cdot}_2^2$ and $\psi_\cY(\cdot)=\frac{1}{2}\norm{\cdot}_2^2$. 
%Note that compactness of $X$ and $Y$ ensures that the sequence ergodic average sequence $\{\tilde x^{(k)},\tilde y^{(k)}\}_{k\geq 1}$ has a limit point and the duality gap defined in \ref{} is a well-defined measure, although perturbation-based variant of the duality gap has been studied in \cite{} which has the advantage of removing such assumption.
%Assumption \ref{assum:bregman}-$(a),(b)$ impose a restriction on the Bregman distance functions. These assumptions have been used in many studies of developing first-order methods, e.g., \cite{hanzely2018accelerated}. Necessary and sufficient conditions in which Assumption \ref{assum:bregman}-$(b)$ holds has been studied in \cite{bauschke2001joint}. For instance, let $\psi_\cX(x)=\sum_{i=1}^m \psi_{x_i}(x_i)$, for $x=[x_i]_{i=1}^m\in\cX=\reals^m$, then for any $x\in X\subseteq \cX$, $\bD_X(x,\cdot)$ is convex if and only if $1/\psi_{x_i}^{''}$ is concave for any $i\in\{1,\hdots,m\}$.  Many examples of Bregman distance used in practice satisfy Assumption \ref{assum:bregman}-$(a),(b)$, e.g. Euclidean distance, squared Mahalanobis distance, generalized Kullback-Leibler divergence, and Itakura-Saito distance.
%\end{remark}

\begin{remark}\label{assum:expected-norm}
Let $(\cU,\norm{\cdot})$ be a finite-dimensional normed real vector space with dual space $\cU^*$. In the analysis of the proposed method we use the following fact: there exists $C_\cU>0$ such that for any vector-valued random variable $W:\Omega\to\cU$, \red{$\mE[\norm{W-\mE[W]}_{\cU^*}^2]\leq  C_\cU\mE[\norm{W}_{\cU^*}^2]$}.

Note that this is a property of the vector space $\cU$ which is true for any finite-dimensional vector space. In more details, 
suppose $\cU$ is a finite-dimensional real vector space equipped with the Euclidean norm denoted by $\norm{\cdot}_2$. Then for any vector-valued random variable $W$, $\mE[\norm{W-\mE[W]}_2^2]=\mE[\norm{W}_2^2]-\mE[\norm{W}_2]^2\leq \mE[\norm{W}_2^2]$; hence, $C_\cU=1$. Moreover, in a finite-dimensional vector space all the norms are equivalent, i.e., for any two arbitrary norms $\norm{\cdot}_\alpha$ and $\norm{\cdot}_\beta$ there exist $c,C>0$ such that $c\norm{\cdot}_\beta\leq\norm{\cdot}_\alpha\leq C\norm{\cdot}_\beta$. Therefore, using the equivalency between an arbitrary norm $\norm{\cdot}_\cU$ and Euclidean norm one can conclude that $C_\cU>0$ exists. It is also worth mentioning that for $(0,+\infty] \ni p$-norms, $C_\cU=m^{\abs{\frac{1}{p}-\frac{1}{q}}}$ where $q\in (0,+\infty]$ is such that $\frac{1}{p}+\frac{1}{q}=1$, and $m$ denotes the dimension of the vector space $\cU$. 
\end{remark}
\begin{assumption}\label{assum:step}
(Step-size conditions) There exist $\alpha,\beta>0$, and $\{\eta^k\}_{k\geq 1}\subset \reals_+$, such that for any $k\geq 1$, the step-sizes $\{\tau^k\}_{k\geq 1}$ and $\{\sigma^k\}_{k\geq 1}$, and the momentum parameters $\{\gamma_x^k,\gamma_y^k,\}_{k\geq 1}\subset (0,1]$ satisfy
\begin{subequations}\label{eq:step}
\begin{align}
&(6C_\cX L_{xx}^2+8C_\cY L_{yx}^2)\eta^k\leq \frac{\gamma_x^k}{\tau^k},&& (6C_\cX L_{xy}^2+8C_\cY L_{yy}^2)\eta^k\leq\frac{\gamma_y^k}{\sigma^k}, \label{eq:step-gamma}\\
& \frac{L^2_{yx}}{\alpha}+2C_\cY L_{yx}^2\eta^k\leq M_x^k, && \frac{L^2_{yy}}{\beta}+2C_\cY L_{yy}^2\eta^k\leq M_y^k, \label{eq:step-tau-sigma}%\\
%& \sqrt{T^k}=\sqrt{T^{k+1}}
\end{align}
\end{subequations}
where $M_x^k\triangleq \frac{1-\gamma_x^k}{\tau^k}- L_{xx}-(6C_\cX L_{xx}^2+8C_\cY L_{yx}^2)\eta^k-\frac{1}{\eta^k}$, and $M_y^k\triangleq \frac{1-\gamma_y^k}{\sigma^k}- (\alpha+\beta)-8C_\cY L_{yy}^2\eta^k-\frac{1}{\eta^k}$. %Moreover, $\{\frac{1}{\eta^k}\}_{k\geq 1}$,  $\{\frac{\gamma_x^k}{\tau^k}T^k,\frac{1-\gamma_x^k}{\tau^k},M_x^k\}_{k\geq 1}$, and $\{\frac{\gamma_y^k}{\sigma^k}T^k, \frac{1-\gamma_y^k}{\sigma^k},M_y^k\}_{k\geq 1}\subset \reals_+$ are non-decreasing sequences.% where $T^k=T(k+1)^2$, for some $T>0$.
\end{assumption}
\begin{remark}\label{rem:step}
Let $L_x\triangleq \sqrt{6C_\cX L_{xx}^2+10C_\cY L_{yx}^2}$ and $L_y\triangleq \sqrt{6C_\cX L_{xy}^2+10C_\cY L_{yy}^2}$. 
The following two choices of the algorithm parameters and the design parameter $\eta^k$ satisfy the step-size conditions in Assumption \ref{assum:step}. 
%A particular choice of the algorithm parameters and the free parameters satisfying Assumption \ref{assum:step} are as follows, 

I) Constant: For $k\geq 1$, $T^k=\bar{T}$, $\tau^k=\tau$, $\sigma^k=\sigma$, $\gamma_x^k=\gamma_x$, $\gamma_y^k=\gamma_y$, and $\eta^k=\eta$, such that
\begin{subequations}\label{eq:specific-step-const}
\begin{align}
&\bar{T} = Tn,\quad \gamma_x=\tfrac{\bar{\gamma}_x}{n},\quad \gamma_y=\tfrac{\bar{\gamma}_y}{n},\quad \alpha=L_{yx},\quad \beta=L_{yy}\\
&\tau=\min\left\{\tfrac{1}{L_x \sqrt{n}},\tfrac{1-\bar{\gamma}_x/n}{L_{xx}+L_{yx}+2L_x}\right\},\quad \sigma=\min\left\{\tfrac{1}{L_y\sqrt{n}},\tfrac{1-\bar{\gamma}_y/n}{L_{yy}+L_{xy}+L_y}\right\}, \label{eq:specific-step-const-ts}\\
&\eta=\min\left\{\tfrac{1}{L_x\sqrt{n}},\tfrac{1}{L_y\sqrt{n}},\tfrac{b_x+\sqrt{b_x^2-4L_x^2}}{2L_x^2},\tfrac{b_y+\sqrt{b_y^2-4L_y^2}}{2L_y^2}\right\},
\end{align}
\end{subequations}
for some $T>0$, $\bar{\gamma}_x,\bar{\gamma}_y\in (0,1)$, where $b_x\triangleq \tfrac{1-1/n}{\tau}-(L_{xx}+L_{yx})$ and $b_y\triangleq \tfrac{1-1/n}{\sigma}-(L_{yy}+L_{xy})$.

II) Non-constant: For any $k\geq 1$,
% \begin{align*}
% &\eta^k=\tfrac{1}{\max\{\sqrt{6C_\cX L_{xx}^2+10C_\cY L_{yx}^2},\sqrt{6C_\cX L_{xy}^2+10C_\cY L_{yy}^2}\}(k+1)}=\Theta\left(\tfrac{1}{k}\right),\\
% &\gamma_x^k=\gamma_y^k=\tfrac{1}{(k+1)^2},\quad \alpha=L_{yx},\quad \beta=L_{yy}\\
% &\tau^k\leq \tfrac{1-(k+1)^{-2}}{L_{xx}+L_{yx}+(6C_\cX L_{xx}^2+10C_\cY L_{yx}^2)\eta^k+1/\eta^k}=\Theta\left(\tfrac{1}{k}\right),\\
% & \sigma^k\leq \tfrac{1-(k+1)^{-2}}{L_{yy}+L_{yx}+10C_\cY L_{yy}^2\eta^k+1/\eta^k}=\Theta\left(\tfrac{1}{k}\right).
% \end{align*}
% \begin{subequations}\label{eq:specific-step}
% \begin{align}
% &\eta^k=\tfrac{1}{{L}(k+1)},\quad \gamma_x^k=\gamma_y^k=\tfrac{1}{(k+1)^2},\quad \alpha=L_{yx},\quad \beta=L_{yy}\\
% &\tau^k=\min\left\{\tfrac{1}{{L}(k+1)},\tfrac{1-(k+1)^{-2}}{2{L}+L_{xx}+L_{yx}}\right\},\quad \sigma^k=\min\left\{\tfrac{1}{{L}(k+1)},\tfrac{1-(k+1)^{-2}}{2{L}+L_{yy}+L_{xy}}\right\}.
% \end{align}
% \end{subequations}
\begin{subequations}\label{eq:specific-step}
\begin{align}
&T^k = T(k+1)^2, \quad \gamma_x^k=\tfrac{\bar{\gamma}_x}{k^2},\quad \gamma_y^k=\tfrac{\bar{\gamma}_y}{k^2},\quad \alpha=L_{yx},\quad \beta=L_{yy}\\
&\tau^k=\min\left\{\tfrac{1}{L_x k},\tfrac{1-\bar{\gamma}_xk^{-2}}{L_{xx}+L_{yx}+2L_x}\right\},\quad \sigma^k=\min\left\{\tfrac{1}{L_yk},\tfrac{1-\bar{\gamma}_yk^{-2}}{L_{yy}+L_{xy}+L_y}\right\},\label{eq:specific-step-ts}\\
&\eta^k=\min\left\{\tfrac{\bar{\gamma}_x}{L_xk},\tfrac{\bar{\gamma}_y}{L_yk},\tfrac{b_x^k+\sqrt{(b_x^k)^2-4L_x^2}}{2L_x^2},\tfrac{b_y^k+\sqrt{(b_y^k)^2-4L_y^2}}{2L_y^2}\right\},
\end{align}
\end{subequations}
for some $T>0$, $\bar{\gamma}_x,\bar{\gamma}_y\in (0,1)$, where $b_x^k\triangleq \frac{1-\gamma_x^k}{\tau^k}-(L_{xx}+L_{yx})$ and $b_y^k\triangleq \frac{1-\gamma_y^k}{\sigma^k}-(L_{yy}+L_{xy})$.

It is worth emphasizing that $\eta^k$ is only a design parameter and need not be computed for running the algorithm.  
\end{remark}

In the following theorem, we state the main result of this paper by providing a bound on the gap function, i.e., $\mathbb{E}[\sup_{z\in\cZ}(\cL(\bar{x},y)-\cL(x,\bar{y}))]$, stating the oracle complexities using constant and non-constant step-sizes in the follow-up corollaries.  

\begin{theorem}\label{thm:main}
Let $\{x_t^k,y_t^k\}_{t,k}$ be  the  sequence  generated  by SVR-APD displayed in Algorithm \ref{alg:SVRPDA} initialized  from  arbitrary  vectors $\tilde{x}^0\in\cX$ and $\tilde{y}^0\in\cY$. Suppose Assumptions \ref{assum:lip} and \ref{assum:bregman} hold, and the step-size sequence $\{\tau^k,\sigma^k\}_{k\geq 1}$ and the momentum parameter sequence $\{\gamma_x^k,\gamma_y^k\}_{k\geq 1}$ satisfy Assumption \ref{assum:step} for some $\alpha,\beta>0$, and $\{\eta^k\}_{k\geq 1}\subset\reals_+$. Moreover, let $z\triangleq (x,y)\in Z\triangleq X\times Y$, then the following holds:

{\bf I)} If the step-sizes and momentum parameters are constant, i.e., $\tau^k=\tau$, $\sigma^k=\sigma$, $\gamma_x^k=\gamma_x$, $\gamma_y^k=\gamma_y$, $T^k=\bar{T}$, then  for $K\geq 1$,  
\begin{align}\label{eq:thm-main-result-constant}
&{\mE\big[\sup_{z\in Z}\big(\cL(\tilde x^{(K)},y)-\cL(x,\tilde y^{(K)})\big)\big]  \leq\tfrac{1}{K\bar{T}}\Big[\Big(\frac{1}{\eta}+\frac{\gamma_x}{\tau}\bar{T}+\frac{1-\gamma_x}{\tau}+M_x\Big)B_X} \nonumber\\ 
& +\Big(\frac{2}{\eta}+\frac{\gamma_y}{\sigma}\bar{T}+\frac{1-\gamma_y}{\sigma}+M_y\Big)B_Y\Big],%\nonumber\\
%&&=\cO\Big(\frac{1}{K^2}\Big),
\end{align}
where $(\tilde{x}^{(K)},\tilde{y}^{(K)})\triangleq \frac{1}{K}\sum_{k=1}^K (\tilde{x}^k,\tilde{y}^k)$. 

{\bf II)} \eyh{If the step-sizes and momentum parameters are non-constant}, then for $K\geq 1$,
\begin{align}\label{eq:thm-main-result}
&{\mE\big[\sup_{z\in Z}\big(\cL(\tilde x^{(K)},y)-\cL(x,\tilde y^{(K)})\big)\big]  \leq\tfrac{1}{S^K}\Big[\Big(\frac{1}{\eta^K}+\frac{\gamma_x^K}{\tau^K}T^K+\frac{1-\gamma_x^K}{\tau^K}+M_x^K\Big)B_X} \nonumber\\ 
& +\Big(\frac{2}{\eta^K}+\frac{\gamma_y^K}{\sigma^K}T^K+\frac{1-\gamma_y^K}{\sigma^K}+M_y^K\Big)B_Y\Big],%\nonumber\\
%&&=\cO\Big(\frac{1}{K^2}\Big),
\end{align}
where $(\tilde{x}^{(K)},\tilde{y}^{(K)})\triangleq \frac{1}{S^K}\sum_{k=1}^K (T^k\tilde{x}^k,T^k\tilde{y}^k)$, $S^K\triangleq \sum_{k=1}^KT^k$. 
%Moreover, if the step-sizes, momentum, and free parameters are selected as in \eqref{eq:specific-step} then $\mE\big[\sup_{z\in Z}\big(\cL(\tilde x^{(K)},y)-\cL(x,\tilde y^{(K)})\big)\big]\leq \cO(1/K^2)$. 
\end{theorem}

\begin{corollary} \label{cor:constant}
Under the premises of Theorem \ref{thm:main} part I, if the step-sizes and parameters of Algorithm \ref{alg:SVRPDA} are selected  as in \eqref{eq:specific-step-const}, then $\mE\big[\sup_{z\in Z}\big(\cL(\tilde x^{(K)},y)-\cL(x,\tilde y^{(K)})\big)\big]  \leq \cO(\frac{1}{K\sqrt{n}})$. Moreover, the oracle complexity is $\cO(\frac{\sqrt{n}}{\epsilon})$.
%the number of iterations required for SVR-APD algorithm to achieve $\epsilon$-gap is $K\geq \cO(\frac{1}{\sqrt{\epsilon}})$, and the number of total primal-dual gradient computations after performing $K$ iterations of SVR-APD algorithm is $nK+\sum_{k=1}^KT^k\leq nK+T(K+1)^3$; hence, the oracle complexity is $\cO(\frac{n}{\sqrt{\epsilon}}+\frac{1}{\epsilon^{1.5}})$.
\end{corollary}

\begin{corollary}\label{cor:nonconstant}
Under the premises of Theorem \ref{thm:main} part II, if the step-sizes and parameters of Algorithm \ref{alg:SVRPDA} are selected as in \eqref{eq:specific-step}, then $\mE\big[\sup_{z\in Z}\big(\cL(\tilde x^{(K)},y)-\cL(x,\tilde y^{(K)})\big)\big]  \leq \cO(\frac{1}{K^2})$. Moreover, the oracle complexity is $\cO(\frac{n}{\sqrt{\epsilon}}+\frac{1}{\epsilon^{1.5}})$.
%the number of iterations required for SVR-APD algorithm to achieve $\epsilon$-gap is $K\geq \cO(\frac{1}{\sqrt{\epsilon}})$, and the number of total primal-dual gradient computations after performing $K$ iterations of SVR-APD algorithm is $nK+\sum_{k=1}^KT^k\leq nK+T(K+1)^3$; hence, the oracle complexity is $\cO(\frac{n}{\sqrt{\epsilon}}+\frac{1}{\epsilon^{1.5}})$.
\end{corollary}

\section{Convergence Analysis}\label{sec:analysis}
To show the convergence results we need the following definitions.

\begin{definition}
We denote expectation and conditional expectation with respect to $\cF^k_t$ by $\mE[\cdot]$ and $\mE[\cdot|\cF^k_t]$, respectively, such that $\cF_t^k\triangleq \sigma(\Psi^k_t)$ where $\sigma(\cdot)$ denotes $\sigma$-algebra and $\Psi^k_t\triangleq \{j^k_0,i^k_0,\hdots,j^k_{t-1},i^k_{t-1}\}$. Similarly we define $\cH_t^k\triangleq \sigma(\Phi^k_t)$ where $\Phi^k_t\triangleq \{j^k_0,i^k_0,\hdots,j^k_{t-1},i^k_{t-1},j^k_t\}$ and let $\mE_{j^k_t}\triangleq \mE[\cdot\mid\cF^k_t]$ and $\mE_{i_t^k}\triangleq \mE[\cdot\mid\cH^k_t]$.
\end{definition}

\begin{definition}\label{def:parameters}
Let $q_t^k\triangleq\grad_y\Phi_{j_k}(x_t^k,y_t^k)-\grad_y\Phi_{j_k}(x_{t-1}^k,y_{t-1}^k)$, $\bar{q}_t^k\triangleq \grad_y\Phi(x_t^k,y_t^k)-\grad_y\Phi(x_{t-1}^k,y_{t-1}^k)$, $\delta_t^y\triangleq \xi_t-\grad_y\Phi(x_t^k,y_t^k)$, and  $\delta^x_t\triangleq\zeta_t-\grad_x\Phi(x_t^k,y_{t+1}^k)$.
\end{definition}

\eyh{\begin{definition}\label{def:mean-bregman}
Let $\mathcal D^k_X(x)\triangleq\tfrac{1}{T^{k-1}}\sum_{\ell=1}^{T^{k-1}}\bD_X(x,x_{\ell}^{k-1})$ and similarly $\mathcal D^k_Y(y)\triangleq\tfrac{1}{T^{k-1}}\sum_{\ell=1}^{T^{k-1}}\bD_Y(y,y_{\ell}^{k-1})$, for $k\geq 1$. 
\end{definition}}

Let $\{T^k\}_{k\geq 1}\subset \integers_+$ be a non-decreasing sequence. We define auxiliary sequences $\{u_t^k\}_{t,k}$, $\{v_t^k\}_{t,k}$, and $\{w_t^k\}_{t,k}$ which are helpful for the analysis of the algorithm. In particular, for any $k\geq 1$ and $t\in \{0,\hdots, T^k-1\}$ we define,
\begin{subequations}
\begin{align}
&u_{t+1}^k \gets \argmin_{x\in X}\{-\fprod{\delta^x_t,x}+\tfrac{1}{\eta^k}\bD_X(x,u_t^k)\}, \label{eq:aux-seq-x}\\
&v_{t+1}^k \gets \argmin_{y\in Y}\{\fprod{\delta_t^y,y}+\tfrac{1}{\eta^k}\bD_Y(y,v_t^k)\},\label{eq:aux-seq-y}\\
&w_{t+1}^k \gets \argmin_{y\in Y}\{\fprod{ q_t^k-\bar{q}_t^k,y}+\tfrac{1}{\eta^k}\bD_Y(y,w_t^k)\},\label{eq:aux-seq-qy}
\end{align}
\end{subequations}
for any $\eta^k>0$, such that $(u_0^0,v_0^0,w_0^0)\triangleq (\tilde{x}^0,\tilde{y}^0,\tilde{y}^0)$ and $(u_{T^k}^k,v_{T^k}^k,w_{T_k}^k)\triangleq (u_0^{k+1},v_0^{k+1},w_0^{k+1})$. 

\begin{lemma}\label{lem:error-inner-prod}
%{\bf Option II}
Let $\{x_t^k,y_t^k\}_{t,k}$ be  the  sequence  generated  by SVR-APD displayed in Algorithm \ref{alg:SVRPDA} initialized  from  arbitrary  vectors $\tilde{x}^0\in\cX$ and $\tilde{y}^0\in\cY$. Let $\{u_t^k,v_t^k,w_t^k\}_{t,k}$ be the auxiliary sequence defined in \eqref{eq:aux-seq-x}-\eqref{eq:aux-seq-qy}. Suppose Assumption \ref{assum:lip} holds and the $\delta_t^k$, $\delta_t^y$, $q_t^k$, and $\bar{q}_t^k$ are defined in Definition \ref{def:parameters}.  For any $x\in X$, $y\in Y$, and $\eta^k>0$ the following results hold for $k\geq 1$ and $t\geq 0$,
\begin{subequations}
\begin{eqnarray}
\lefteqn{\fprod{\delta_t^x,x-x_{t+1}^k}\leq} \\ \nonumber 
&&\frac{1}{\eta^k}\Big(\bD_X(x,u_t^k)-\bD_X(x,u_{t+1}^k)\Big)+\eta^k\norm{\delta_t^x}_{\cX^*}^2+\fprod{\delta_t^x,u_t^k-x_t^k}+\frac{1}{\eta^k}\bD_X(x_{t+1}^k,x_t^k) , \label{eq:bound-delta-x}\\
\lefteqn{\fprod{\delta_t^y,y_{t+1}^k-y}\leq}  \\ \nonumber
&&\frac{1}{\eta^k}\Big(\bD_Y(y,v_t^k)-\bD_Y(y,v_{t+1}^k)\Big)+\eta^k\norm{\delta_t^y}_{\cY^*}^2+\fprod{\delta_t^y,y_t^k-v_t^k}+\frac{1}{\eta^k}\bD(y_{t+1}^k,y_t^k), \label{eq:bound-delta-y}\\
\lefteqn{\fprod{q_t^k-\bar{q}_t^k,y_t^k-y}\leq }\\
\nonumber&&\frac{1}{\eta^k}\Big(\bD_Y(y,w_t^k)-\bD_Y(y,w_{t+1}^k)\Big)+\frac{\eta^k}{2}\norm{q_t^k-\bar{q}_t^k}_{\cY^*}^2+\fprod{q_t^k-\bar{q}_t^k,y_t^k-w_t^k}. \label{eq:bound-qy}
\end{eqnarray}
\end{subequations}
%\red{I change $\eta^k^k$ to $1/\eta^k^k$ in above, change next two equations accordingly}
\end{lemma}
{\it Proof} 
We split the inner product into $\fprod{\delta_t^x,x-x_t^k}$ and $\fprod{\delta^x_t,x_t^k-x_{t+1}^k}$ and provide an upper bound for each term. Using Lemma \ref{lem:triangle-breg}-$(c)$ for \eqref{eq:aux-seq-x} with $s=\delta_t^x$, $t=1/\eta^k$, and $f\equiv 0$, we conclude that for any $x\in X$,
\begin{align}\label{eq:bound-xu}
{\fprod{\delta_t^x,x-x_t^k}} = &\fprod{\delta_t^x,x-u_t^k}+\fprod{\delta_t^x,u_t^k-x_t^k}\\
\leq & \frac{1}{\eta^k}\Big(\bD_X(x,u_t^k)-\bD_X(x,u_{t+1}^k)\Big)+\frac{\eta^k}{2}\norm{\delta_t^x}_{\cX^*}^2+\fprod{\delta_t^x,u_t^k-x_t^k}. \nonumber
\end{align}
Moreover, using Young's inequality and strong convexity of the Bregman distance function we have  
\begin{align}\label{eq:inner-prod-error}
\fprod{\delta^x_t,x_t^k-x_{t+1}^k}\leq \frac{\eta^k}{2}\norm{\delta_t^x}_{\cX^*}^2+\frac{1}{\eta^k}\bD_X(x_{t+1}^k,x_t^k) 
\end{align}
Adding \eqref{eq:inner-prod-error} to \eqref{eq:bound-xu} gives \eqref{eq:bound-delta-x}. Similarly, using \eqref{eq:aux-seq-y} and \eqref{eq:aux-seq-qy} one can obtain the results in  \eqref{eq:bound-delta-y} and \eqref{eq:bound-qy}, respectively.\qed

In the following lemma, we derive upper bounds for the error of estimating gradients.

\begin{lemma}\label{lem:upper-noise}
Under the premises of Lemma \ref{lem:error-inner-prod}, for any $k\geq 1$ and $t\geq 0$, $\mE_{i_t^k}[\delta^x_t]=\mE_{j_t^k}[\delta^y_t]=\mE_{j_t^k}[q_t^k-\bar{q}_t^k]=0$, the following hold for some $C_\cX,C_\cY>0$.
\begin{subequations}\label{eq:norm-error}
\begin{align}
%&\mE_{i_t^k}\Big[\norm{\delta^x_t}_{\cX^*}^2\Big] \label{eq:norm-error-x}  \leq  \red{6C_\cX}\mE_{i_t^k}\Big[ L_{xx}^2\bD_X(x_{t+1}^k,x_t^k)+ L_{xx}^2\bD_X(x_{t+1}^k,\tilde x^{k-1})\nonumber\\
%&+ L_{xy}^2\bD_Y(y_{t+1}^k,\tilde y^{k-1})\Big],  \\
&\eyh{\mE_{i_t^k}\Big[\norm{\delta^x_t}_{\cX^*}^2\Big] \label{eq:norm-error-x}  \leq  \red{6C_\cX}\mE_{i_t^k}\Big[ L_{xx}^2\bD_X(x_{t+1}^k,x_t^k)+ L_{xx}^2\cD_X^k(x_{t+1}^k)+ L_{xy}^2\cD_Y^k(y_{t+1}^k)\Big],}  \\
&\mE_{j_t^k}\Big[\norm{\delta^y_t}_{\cY^*}^2\Big]  \label{eq:norm-error-y}\leq \red{8C_\cY}\mE_{j_t^k}\Big[L_{yx}^2\bD_X(x_{t+1}^k,x_t^k)+L_{yy}^2\bD_Y(y_{t+1}^k,y_t^k)\nonumber\\
%&+L_{yy}^2\bD_Y(y_{t+1}^k,\tilde y^{k-1})+L_{yx}^2\bD_X(x_{t+1}^k,\tilde x^{k-1})\Big], \\
&\quad \eyh{+L_{yy}^2\cD_Y^k(y_{t+1}^k)+L_{yx}^2\cD_X^k(x_{t+1}^k)\Big],}\\
&\mE_{j_t^k}\Big[\norm{q_t^k-\bar{q}_t^k}_{\cY^*}^2\Big]   \leq  \red{4C_\cY}\mE_{j_t^k}\Big[L_{yx}^2\bD_X(x_{t}^k,x_{t-1}^k)+L_{yy}^2\bD_Y(y_{t}^k,y_{t-1}^k)\Big],  \label{eq:norm-error-q}
\end{align}
\end{subequations}
% \begin{eqnarray}
% \lefteqn{\mE_{i_t^k}\Big[\fprod{\zeta_t-\grad_x\Phi(x_t^k,y_{t+1}^k),x-x_{t+1}^k}\Big]\leq}\\ 
% &&\mE_{i_t^k}\Big[\frac{2(L_{xx}\eta^k^k)^2+1}{2\eta^k^k}\norm{x_t^k-x_{t+1}^k}^2+{2 L_{xx}^2}{\eta^k^k}\norm{x_{t+1}^k-\tilde x^{k-1}}^2+{2 L_{xy}^2}{\eta^k^k}\norm{y_{t+1}^k-\tilde y^{k-1}}^2\Big] \nonumber \\
% && +\mE_{i_t^k}\Big[\frac{1}{2\eta^k^k}\Big(\norm{x-u_t^k}^2-\norm{x-u_{t+1}^k}^2\Big) +\fprod{\delta_t^x,u_t^k-x_t^k}\Big]
% \end{eqnarray}
% \begin{eqnarray}
% \lefteqn{\mE_{j_t^k}\Big[\fprod{\xi_t-\grad_y\Phi(x_t^k,y_t^k),y_{t+1}^k-y}\Big]\leq}\\ 
% &&\mE_{j_t^k}\Big[\frac{4L_{yy}^2+\eta^k^2}{2\eta^k}\norm{y_t^k-y_{t+1}^k}^2+\frac{2 L_{yy}^2}{\eta^k}\norm{x_t^k-x_{t+1}^k}^2\nonumber\\
% &&+\frac{2 L_{yy}^2}{\eta^k}\norm{y_{t+1}^k-\tilde y^{k-1}}^2+\frac{2 L_{yx}^2}{\eta^k}\norm{x_{t+1}^k-\tilde x^{k-1}}^2\Big] \nonumber
% \end{eqnarray}
\end{lemma}
{\it Proof} 
The unbiasedness of the stochastic noises $\delta_t^x$, $\delta^y_t$, and $q_t^k-\bar{q}_t^k$ clearly holds due to the uniform sampling of the sum-function components. Next,
from the definition of ${\delta_t^x}$ we have that
\begin{align}
\mE_{i_t^k}\Big[\norm{\delta^x_t}_{\cX^*}^2\Big]  &  =\mE_{i_t^k}\Big[\Big\|(\grad_x\Phi_{i_t^k}(x_t^k,y_{t+1}^k)-\grad_x\Phi_{i_t^k}(\tilde x^{k-1},\tilde y^{k-1}))\nonumber\\
&-(\grad_x\Phi(x_t^k,y_{t+1}^k)-\grad_x\Phi(\tilde x^{k-1},\tilde y^{k-1}))\Big\|_{\cX^*}^2\Big]\nonumber \\
& \leq \red{C_\cX}\mE_{i_t^k}\Big[\norm{\grad_x\Phi_{i_t^k}(x_t^k,y_{t+1}^k)-\grad_x\Phi_{i_t^k}(\tilde{x}^{k-1},\tilde y^{k-1})}_{\cX^*}^2\Big],\label{eq:bound-delta-1}
\end{align}
where the inequality is concluded using Remark \ref{assum:expected-norm} for $W=\grad_x\Phi_{i_t^k}(x_t^k,y_{t+1}^k)-\grad_x\Phi_{i_t^k}(\tilde{x}^{k-1},\tilde y^{k-1})$.
Next, using the triangle inequality, the fact that for any $\{a_i\}_{i=1}^m\subset\reals$, $(\sum_{i=1}^m a_i)^2\leq m\sum_{i=1}^ma_i^2$, and Lipschitz continuity of $\grad_x \Phi_i$ in Assumption \ref{assum:lip} for $i\in\cN$ we conclude that
\begin{eqnarray}\label{eq:bound-delta-2}
\lefteqn{\mE_{i_t^k}\Big[\norm{\grad_x\Phi_{i_t^k}(x_t^k,y_{t+1}^k)-\grad_x\Phi_{i_t^k}(\tilde{x}^{k-1},\tilde y^{k-1})}_{\cX^*}^2\Big]}  \\
&& \leq 3\mE_{i_t^k}\Big[\norm{\grad_x\Phi_{i_t^k}(x_t^k,y_{t+1}^k)-\grad_x\Phi_{i_t^k}(x_{t+1}^k,y_{t+1}^k)}_{\cX^*}^2\nonumber\\&&
\quad+\norm{\grad_x\Phi_{i_t^k}(x_{t+1}^k,y_{t+1}^k)-\grad_x\Phi_{i_t^k}(x_{t+1}^k,\tilde y^{k-1})}_{\cX^*}^2\nonumber \\
&&\quad+\norm{\grad_x\Phi_{i_t^k}(x_{t+1}^k,\tilde y^{k-1})-\grad_x\Phi_{i_t^k}(\tilde{x}^{k-1},\tilde y^{k-1})}_{\cX^*}^2\Big]\nonumber \\ \nonumber
&&\leq 3\mE_{i_t^k}\Big[ L_{xx}^2\norm{x_t^k-x_{t+1}^k}_{\cX}^2+ L_{xx}^2\norm{x_{t+1}^k-\tilde x^{k-1}}_{\cX}^2+ L_{xy}^2\norm{y_{t+1}^k-\tilde y^{k-1}}_{\cY}^2\Big]\nonumber\\
&&\leq \eyh{3\mE_{i_t^k}\Big[ L_{xx}^2\norm{x_t^k-x_{t+1}^k}_{\cX}^2+\tfrac{1}{T^{k-1}}\sum_{\ell=1}^{T^{k-1}}\big( L_{xx}^2\norm{x_{t+1}^k- x_{\ell}^{k-1}}_{\cX}^2}\nonumber\\
&&\quad \eyh{+ L_{xy}^2\norm{y_{t+1}^k- y_{\ell}^{k-1}}_{\cY}^2\big)\Big]}\nonumber.
\end{eqnarray}
Finally, combining \eqref{eq:bound-delta-2} with \eqref{eq:bound-delta-1}, using strong convexity of the Bregman distance function, \eyh{and Definition of \eqref{def:mean-bregman}} lead to \eqref{eq:norm-error-x}. 

For proving the inequalities in \eqref{eq:norm-error-y} and \eqref{eq:norm-error-q} one can use a similar argument. 
In particular, 
\begin{align}
 &{\mE_{j_t^k}\Big[\norm{\delta_t^y}_{\cY^*}^2\Big]}\nonumber\\
 &= \mE_{j_t^k}\Big[\Big\|\grad_y\Phi_{j_t^k}(x_t^k,y_t^k)-\grad_y\Phi_{j_t^k}(\tilde{x}^{k-1},\tilde{y}^{k-1})-(\grad_y\Phi(x_t^k,y_t^k)-\grad_y\Phi(\tilde{x}^{k-1},\tilde{y}^{k-1}))\Big\|_{\cY^*}^2 \Big]\nonumber\\
 &\leq  \red{C_\cY}\mE_{j_t^k}\Big[\norm{\grad_y\Phi_{j_t^k}(x_t^k,y_t^k)-\grad_y\Phi_{j_t^k}(\tilde{x}^{k-1},\tilde{y}^{k-1})}^2_{\cY^*}\Big]\nonumber\\
%  &\leq \red{4C_\cY}\mE_{j_t^k}\Big[\norm{\grad_y\Phi_{j_t^k}(x_t^k,y_t^k)-\grad_y\Phi_{j_t^k}(x_{t+1}^k,y_t^k)}^2_{\cY^*}\nonumber\\
%  &\quad +\norm{\grad_y\Phi_{j_t^k}(x_{t+1}^k,y_t^k)-\grad_y\Phi_{j_t^k}(x_{t+1}^k,y_{t+1}^k)}^2_{\cY^*} \nonumber\\
%  &\quad +\norm{\grad_y\Phi_{j_t^k}(x_{t+1}^k,y_{t+1}^k)-\grad_y\Phi_{j_t^k}(x_{t+1}^k,\tilde{y}^{k-1})}^2_{\cY^*}\nonumber\\
%  &\quad+\norm{\grad_y\Phi_{j_t^k}(x_{t+1}^k,\tilde{y}^{k-1})-\grad_y\Phi_{j_t^k}(\tilde{x}^{k-1},\tilde{y}^{k-1})}^2_{\cY^*}\Big]\nonumber\\
 &\leq \red{4C_\cY}\mE_{j_t^k}\Big[ L_{yx}^2\norm{x_t^k-x_{t+1}^k}_{\cX}^2+L_{yy}^2\norm{y_t^k-y_{t+1}^k}_\cY^2+ L_{yy}^2\norm{y_{t+1}^k-\tilde y^{k-1}}_{\cY}^2\nonumber\\
 &\quad+L_{yx}^2\norm{x_{t+1}^k-\tilde x^{k-1}}_{\cX}^2\Big]\nonumber\\
 &\leq \eyh{\red{4C_\cY}\mE_{j_t^k}\Big[ L_{yx}^2\norm{x_t^k-x_{t+1}^k}_{\cX}^2+L_{yy}^2\norm{y_t^k-y_{t+1}^k}_\cY^2}\nonumber\\
 &\quad\eyh{+ \tfrac{1}{T^{k-1}}\sum_{\ell=1}^{T^{k-1}}\big(L_{yy}^2\norm{y_{t+1}^k-y_\ell^{k-1}}_{\cY}^2+L_{yx}^2\norm{x_{t+1}^k-x_\ell^{k-1}}_{\cX}^2\Big]}\nonumber.
\end{align}
Moreover, 
\begin{align}
 {\mE_{j_t^k}\Big[\norm{q_t^k-\bar{q}_t^k}_{\cY^*}^2\Big]}
 & = \mE_{j_t^k}\Big[\Big\|\grad_y\Phi_{j_t^k}(x_t^k,y_t^k)-\grad_y\Phi_{j_t^k}(x_{t-1}^k,y_{t-1}^k)\nonumber\\&-(\grad_y\Phi(x_t^k,y_t^k)-\grad_y\Phi(x_{t-1}^k,y_{t-1}^k)\Big\|_{\cY^*}^2 \Big]\nonumber\\
 &\leq \red{C_\cY}\mE_{j_t^k}\Big[\norm{\grad_y\Phi_{j_t^k}(x_t^k,y_t^k)-\grad_y\Phi_{j_t^k}(x_{t-1}^k,y_{t-1}^k)}_{\cY^*}^2\Big]\nonumber\\
 &\leq \red{2C_\cY}\mE_{j_t^k}\Big[\norm{x_t^k-x_{t-1}^k}_\cX^2+\norm{y_t^k-y_{t-1}^k}^2_\cY\Big].\nonumber  \qquad\qquad \qed 
\end{align}
 %Finally, combining \eqref{eq:bound-xu}, \eqref{eq:inner-prod-error}, and \eqref{eq:norm-error} implies the desired result in \eqref{}. Moreover, a similar analysis can be used to show \eqref{}. 

\begin{lemma}\label{lem:inner-prod}
Suppose Assumption \ref{assum:lip} holds, then the following inequality holds for any $x\in X,\bar{x}\in X^\circ$, $y,\tilde{y}\in Y$, and $\bar{y}\in Y^\circ$,
\begin{align}\label{eq:inner-prod}
\abs{\fprod{\grad_y\Phi(x,y)-\grad_y\Phi(\bar{x},\bar{y}),y-\tilde y}}&
\leq \frac{L_{yx}^2}{\alpha}\bD_X(x,\bar{x})+\frac{L_{yy}^2}{\beta}\bD_Y(y,\bar{y})\nonumber\\&+(\alpha+\beta)\bD_Y(\tilde y,y),
\end{align}
for any $\alpha,\beta>0$. Moreover, if $L_{yy}=0$, then for any $\alpha>0$:
\begin{align}\label{eq:inner-prod-beta-0}
\abs{\fprod{\grad_y\Phi(x,y)-\grad_y\Phi(\bar{x},\bar{y}),y-\tilde y}}
\leq \frac{L_{yx}^2}{\alpha}\bD_X(x,\bar{x})+\alpha\bD_X(y,\tilde y),
\end{align}

\end{lemma}
{\it Proof}  By adding and subtracting $\grad_y\Phi(\bar{x},y)$ and triangle inequality we obtain 
\begin{eqnarray}\label{eq:tri-ineq}
\lefteqn{\abs{\fprod{\grad_y\Phi(x,y)-\grad_y\Phi(\bar{x},\bar{y}),y-\tilde y}}}\nonumber \\
&&\leq \abs{\fprod{\grad_y\Phi(x,y)-\grad_y\Phi(\bar{x},y),y-\tilde y}}+\abs{\fprod{\grad_y\Phi(\bar{x},y)-\grad_y\Phi(\bar{x},\bar{y}),y-\tilde y}}.%\\
%&&\leq \frac{L_{yx}^2}{2\alpha}\norm{x-\bar{x}}_\cX^2+\frac{\alpha}{2}\norm{y-\tilde y}_\cY^2+\frac{L_{yy}^2}{2\beta}\norm{y-\bar{y}}_\cY^2+\frac{\beta}{2}\norm{y-\tilde y}_\cY^2
\end{eqnarray}
Next, using Young's inequality and Assumption \ref{assum:lip} we conclude that 
\begin{eqnarray*}
\lefteqn{\abs{\fprod{\grad_y\Phi(x,y)-\grad_y\Phi(\bar{x},\bar{y}),y-\tilde y}}}\\
&&\leq \frac{L_{yx}^2}{2\alpha}\norm{x-\bar{x}}_\cX^2+\frac{\alpha}{2}\norm{y-\tilde y}_\cY^2+\frac{L_{yy}^2}{2\beta}\norm{y-\bar{y}}_\cY^2+\frac{\beta}{2}\norm{y-\tilde y}_\cY^2.
\end{eqnarray*}
The result in \eqref{eq:inner-prod} immediately follows using strong convexity of Bregman distance functions. Moreover, if $L_{yy}=0$, then the second inner product in the right hand side of \eqref{eq:tri-ineq} will be zero and the result can be concluded. \qed

An immediate consequence of Lemma \ref{lem:inner-prod} by setting $(x,y)=(x_t^k,y_t^k)$ and $(\bar{x},\bar{y})=(x_{t-1}^k,y_{t-1}^k)$ is that for any $\tilde{y}\in\cY$, $k\geq 1$, and $t\geq 0$
\begin{align}\label{eq:inner-prod-sum}
\abs{\fprod{\bar{q}_t^k,y_t^k-\tilde y}}
\leq \frac{L_{yx}^2}{\alpha}\bD_X(x_t^k,x_{t-1}^k)+\frac{L_{yy}^2}{\beta}\bD_Y(y_t^k,y_{t-1}^k)+(\alpha+\beta)\bD_Y(\tilde y,y_t^k).
\end{align}

% Important property:
% \begin{equation}\label{eq:important-eq}
% \fprod{a-b,c-d}=\frac{1}{2}\Big(\norm{a-d}^2+\norm{b-c}^2-\norm{a-c}^2-\norm{b-d}^2\Big)
% \end{equation}
In the next lemma, we provide a one-step analysis for SVR-APD which is the main building block for showing the rate result stated in Theorem \ref{thm:main}.

\begin{lemma}\label{lem:one-step}
Let $\{x_t^k,y_t^k\}_{t,k}$ be  the  sequence  generated  by SVR-APD displayed in Algorithm \ref{alg:SVRPDA} initialized  from  arbitrary  vectors $\tilde{x}^0\in\cX$ and $\tilde{y}^0\in\cY$. Let $\{u_t^k,v_t^k,w_t^k\}_{t,k}$ be the auxiliary sequence defined in \eqref{eq:aux-seq-x}-\eqref{eq:aux-seq-qy}. Suppose Assumption \ref{assum:lip} holds and the $\delta_t^k$, $\delta_t^y$, $q_t^k$, and $\bar{q}_t^k$ are defined in Definition \ref{def:parameters}. Then for any $x\in X$, $y\in Y$, and $\eta^k>0$ the following inequality holds for $k\geq 1$ and $t\in \{0,\hdots, T^k-1\}$,
\begin{align}\label{eq:one-step}
&\cL(x_{t+1}^k,y)-\cL(x,y_{t+1}^k)\leq \\
&A_t^k-B_t^k +\fprod{\bar{q}_t^k,y_t^k-y}-\fprod{\bar{q}^k_{t+1},y_{t+1}^k-y}+\frac{1}{\eta^k}\Big(\bD_X(x,u_t^k)-\bD_X(x,u_{t+1}^k)\Big) \nonumber \\
& +\frac{1}{\eta^k}\Big(\bD_Y(y,v_t^k)-\bD_Y(y,v_{t+1}^k)\Big) +\frac{1}{\eta^k}\Big(\bD_Y(y,w_t^k)-\bD_Y(y,w_{t+1}^k)\Big) \nonumber \\
&+\frac{\gamma_x^k}{\tau^k}\Big(\eyh{\cD_X^k(x)} -\bD_X(x,x_{t+1}^k)\Big)+\frac{1-\gamma_x^k}{\tau^k}\Big(\bD_X(x,x_t^k)-\bD_X(x,x_{t+1}^k)\Big) \nonumber\\
&+\big((6C_\cX L_{xx}^2+8C_\cY L_{yx}^2)\eta^k-\frac{\gamma_x^k}{\tau^k}\big)\eyh{\cD_X^k(x_{t+1}^k)} \nonumber\\
&+(\frac{L_{yx}^2}{\alpha}+2C_\cY L_{yx}^2\eta^k) \bD_X(x_t^k,x_{t-1}^k)-{M_x^k}\bD_X(x_{t+1}^k,x_t^k) \nonumber\\
& +\frac{\gamma_y^k}{\sigma^k}\Big(\eyh{\cD_Y^k(y)} -\bD_Y(y,y_{t+1}^k)\Big)+\frac{1-\gamma_y^k}{\sigma^k}\Big(\bD_Y(y,y_t^k)-\bD_Y(y,y_{t+1}^k)\Big)\nonumber\\
&+\big((6C_\cX L_{xy}^2+8C_\cY L_{yy}^2)\eta^k-\frac{\gamma_y^k}{\sigma^k}\big)\eyh{\cD_Y^k(y_{t+1}^k)}\nonumber \\
&+(\frac{L_{yy}^2}{\beta}+2C_\cY L_{yy}^2\eta^k) \bD_Y(y_t^k,y_{t-1}^k)-{M_y^k}\bD_Y(y_{t+1}^k,y_t^k), \nonumber\\
&{A_t^k\triangleq} \fprod{\delta_t^x,u_t^k-x_t^k}+\fprod{\delta_t^y,y_t^k-v_t^k}+\fprod{q_t^k-\bar{q}_t^k,y_t^k-w_t^k}\\
&\qquad+\eta^k\big(\norm{\delta_t^x}_{\cX^*}^2+\norm{\delta_t^y}_{\cY^*}^2 +\tfrac{1}{2}\norm{q_t^k-\bar{q}_t^k}_{\cY^*}^2\big), \nonumber\\
&B_t^k\triangleq \eta^k\Big((6C_\cX L_{xx}^2+8C_\cY L_{yx}^2)(\bD_X(x_{t+1}^k,x_t^k)+\eyh{\cD_X^k(x_{t+1}^k)})\\
&\qquad+8C_\cY L_{yy}^2\bD_Y(y_{t+1}^k,y_t^k)+(6C_\cX L_{xy}^2+8C_\cY L_{yy}^2)\eyh{\cD_Y^k(y_{t+1}^k)}\nonumber\\
&\qquad+2C_\cY (L_{yx}^2\bD_X(x_t^k,x_{t-1}^k)+L_{yy}^2\bD_Y(y_t^k,y_{t-1}^k))\Big), \nonumber
\end{align}
where $M_x^k$ and $M_y^k$ are defined in Assumption \ref{assum:step}.
%$\bar{q}_t^k\triangleq \grad_y\Phi(x_t^k,y_t^k)-\grad_y\Phi(x_{t-1}^k,y_{t-1}^k)$, $\delta_t^y\triangleq \xi_t-\grad_y\Phi(x_t^k,y_t^k)$, $\delta^x_t\triangleq\zeta_t-\grad_x\Phi(x_t^k,y_{t+1}^k)$,
%$M_x^k\triangleq \frac{1-\gamma_x^k}{\tau^k}- L_{xx}-(6C_\cX L_{xx}^2+8C_\cY L_{yx}^2)\eta^k-\frac{1}{\eta^k}$, and $M_y^k\triangleq \frac{1-\gamma_y^k}{\sigma^k}- (\alpha+\beta)-8C_\cY L_{yy}^2\eta^k-\frac{1}{\eta^k}$.
\end{lemma}
{\it Proof} 
Applying Lemma \ref{lem:triangle-breg}-$(a)$, on the update rule of $y_{t+1}^k$, implies that for any $y\in Y$,
\begin{align*}
h(y_{t+1}^k)-h(y)&\leq 
\fprod{\xi_t+  q_t^k,y_{t+1}^k-y}
%+\frac{\gamma_y^k}{\sigma^k}\fprod{y_t^k-\tilde y^{k-1},y-y_{t+1}^k} \label{eq:prox-h}\\
%&+\frac{1}{\sigma^k}\Big(\bD_Y(y,y_t^k)-\bD_Y(y,y_{t+1}^k)-\bD_Y(y_{t+1}^k,y_t^k)\Big) \nonumber
+\frac{1}{\sigma^k}\fprod{\grad\psi_\cY(y_{t+1}^k)-\grad\psi_\cY(\hat{y}_t^k),y-y_{t+1}^k}\nonumber\\
&= \fprod{\xi_t+  q_t^k,y_{t+1}^k-y} \\
&\quad + \frac{1}{\sigma^k}\fprod{\grad\psi_\cY(y_{t+1}^k)-(1-\gamma_y^k)\grad\psi_\cY({y}_t^k)-\gamma_y^k\eyh{s^{k-1}},y-y_{t+1}^k}\\
&=\fprod{\xi_t+  q_t^k,y_{t+1}^k-y} +
\frac{1-\gamma_y^k}{\sigma^k}\fprod{\grad\psi_\cY(y_{t+1}^k)-\grad\psi_\cY(y_t^k),y-y_{t+1}^k}\\
&\quad \eyh{+ \frac{\gamma_y^k}{\sigma^kT^{k-1}}\sum_{\ell=1}^{T^{k-1}}\fprod{\grad\psi_\cY(y_{t+1}^k)-\grad\psi_\cY(y_{\ell}^{k-1}),y-y_{t+1}^k},}
\end{align*}
where in the first equality we used Lemma \ref{lem:triangle-breg}-$(d)$ and update rule of $\hat{y}^k_t$ in line 7 of Algorithm \ref{alg:SVRPDA}. Using the generalized three-point property of Bregman distance in Lemma \ref{lem:triangle-breg}-$(b)$ twice; for $\fprod{\grad\psi_\cY(y_{t+1}^k)-\grad\psi_\cY(y_t^k),y-y_{t+1}^k}$ and $\fprod{\grad\psi_\cY(y_{t+1}^k)-\eyh{\grad\psi_\cY(y_{\ell}^{k-1})},y-y_{t+1}^k}$, the last inequality can be written as
\begin{align}
h(y_{t+1}^k)-h(y)\leq 
& \fprod{\xi_t+ q_t^k,y_{t+1}^k-y} \nonumber\\
&+ \frac{1-\gamma_y^k}{\sigma^k}\Big(\bD_Y(y,y_t^k)-\bD_Y(y,y_{t+1}^k)-\bD_Y(y_{t+1}^k,y_t^k)\Big) \nonumber\\
&+\eyh{\frac{\gamma_y^k}{\sigma^kT^{k-1}}\sum_{\ell=1}^{T^{k-1}}\Big(\bD_Y(y,y_{\ell}^{k-1})-\bD_Y(y,y_{t+1}^k)-\bD_Y(y_{t+1}^k,y_{\ell}^{k-1})\Big)}.\label{eq:prox-h}
\end{align}
Using the definition of $\delta_t^y$, i.e., $\delta_t^y= \xi_t-\grad_y\Phi(x_t^k,y_t^k)$, Definition \ref{def:mean-bregman}, and rearranging the terms we obtain that 
\begin{align}\label{eq:prox-h-acc}
h(y_{t+1}^k)-h(y)\leq & \fprod{\delta_t^y,y_{t+1}^k-y}+\fprod{\grad_y\Phi(x_t^k,y_t^k),y_{t+1}^k-y} \nonumber\\
&+\fprod{ q_t^k,y_t^k-y}+\fprod{ q_t^k,y_{t+1}^k-y_t^k} \nonumber \\
&+ \frac{1-\gamma_y^k}{\sigma^k}\Big(\bD_Y(y,y_t^k)-\bD_Y(y,y_{t+1}^k)-\bD_Y(y_{t+1}^k,y_t^k)\Big) \nonumber\\
&+\eyh{\frac{\gamma_y^k}{\sigma^k}\Big(\cD_Y^k(y)-\bD_Y(y,y_{t+1}^k)-\cD_Y^k(y_{t+1}^k)\Big)} .
\end{align}
Adding $\fprod{\grad_y\Phi(x_{t+1}^k,y_{t+1}^k),y-y_{t+1}^k}$ to both sides of \eqref{eq:prox-h-acc}  and adding and subtracting $\fprod{\bar{q}_t^k,y_t^k-y}$ to the right hand side of \eqref{eq:prox-h-acc} lead to
\begin{eqnarray}\label{eq:prox-h-acc-2}
\lefteqn{h(y_{t+1}^k)-h(y)+\fprod{\grad_y\Phi(x_{t+1}^k,y_{t+1}^k),y-y_{t+1}^k}\leq}\nonumber\\
&&\fprod{\delta_t^y,y_{t+1}^k-y}-\fprod{\bar{q}_{t+1}^k,y_{t+1}^k-y}\nonumber \\
&&+\big(\fprod{q_t^k-\bar{q}_t^k,y_t^k-y}+ \fprod{\bar{q}_t^k,y_t^k-y}+\fprod{q_t^k,y_{t+1}^k-y_t^k}\big) \nonumber \\
&&+ \frac{1-\gamma_y^k}{\sigma^k}\Big(\bD_Y(y,y_t^k)-\bD_Y(y,y_{t+1}^k)-\bD_Y(y_{t+1}^k,y_t^k)\Big) \nonumber\\
&&+\eyh{\frac{\gamma_y^k}{\sigma^k}\Big(\cD_Y^k(y)-\bD_Y(y,y_{t+1}^k)-\cD_Y^k(y_{t+1}^k)\Big)} .
\end{eqnarray}
Now using \eqref{eq:inner-prod-sum} we obtain 
\begin{align}\label{eq:inner-prod-q}
\fprod{q_t^k,y_{t+1}^k-y_t^k}\leq \frac{L_{yx}^2}{\alpha}\bD_X(x_t^k,x_{t-1}^k)+\frac{L_{yy}^2}{\beta}\bD_Y(y_t^k,y_{t-1}^k)+(\alpha+\beta)\bD_Y(y_{t+1}^k,y_t^k) ,
\end{align}
for any $\alpha,\beta>0$, and note that if $L_{yy}=0$, then $\beta=0$ and we define $0^2/0=0$.
Moreover, $\Phi(x,\cdot)$ is a concave function, for any $x\in\cX$, therefore, we have that 
\begin{align}\label{eq:concave-phi}
\Phi(x_{t+1}^k,y)-\Phi(x_{t+1}^k,y_{t+1}^k)\leq \fprod{\grad_y\Phi(x_{t+1}^k,y_{t+1}^k),y-y_{t+1}^k}.
\end{align} 
Using \eqref{eq:inner-prod-q} within \eqref{eq:prox-h-acc-2} and adding \eqref{eq:concave-phi} to the result implies that
\begin{eqnarray}\label{eq:prox-h-acc-3}
\lefteqn{h(y_{t+1}^k)-h(y)+\Phi(x_{t+1}^k,y)-\Phi(x_{t+1}^k,y_{t+1}^k)\leq}\\
&& \fprod{\delta_t^y,y_{t+1}^k-y}+\fprod{q_t^k-\bar{q}_t^k,y_t^k-y}+ \fprod{\bar{q}_t^k,y_t^k-y}-\fprod{\bar{q}_{t+1}^k,y_{t+1}^k-y} \nonumber \\
&&+\frac{1-\gamma_y^k}{\sigma^k}\Big(\bD_Y(y,y_t^k)-\bD_Y(y,y_{t+1}^k)-\bD_Y(y_{t+1}^k,y_t^k)\Big) \nonumber\\
&&+\eyh{\frac{\gamma_y^k}{\sigma^k}\Big(\cD_Y^k(y)-\bD_Y(y,y_{t+1}^k)-\cD_Y^k(y_{t+1}^k)\Big)}  \nonumber \\
&&+ \Big(\frac{L_{yx}^2}{\alpha}\bD_X(x_t^k,x_{t-1}^k)+\frac{L_{yy}^2}{\beta}\bD_Y(y_t^k,y_{t-1}^k)+(\alpha+\beta)\bD_Y(y_{t+1}^k,y_t^k)\Big). \nonumber
\end{eqnarray}

Next, we analyze the update rule of $x_{t+1}^k$. Applying Lemma \ref{lem:triangle-breg}-$(a)$  on the update rule of $x_{t+1}^k$, implies that for any $x\in X$,
\begin{align}
f(x_{t+1}^k)-f(x)\leq &\fprod{\zeta_t,x-x_{t+1}^k}
+\frac{1}{\tau^k}\fprod{\grad\psi_\cX(x_{t+1}^k)-\grad\psi_\cX(\hat{x}^k_t),x-x_{t+1}^k}\nonumber \\
=&  \fprod{\zeta_t,x-x_{t+1}^k}\eyh{+\frac{1-\gamma_x^k}{\tau^k}\fprod{\grad\psi_\cX(x_{t+1}^k)-\grad\psi_\cX(x_t^k),x-x_{t+1}^k}}\nonumber\\
&
\eyh{+\frac{\gamma_x^k}{\tau^k}\fprod{\grad\psi_\cX({x}_{t+1}^k)-r^{k-1},x-x_{t+1}^k},}
%+\frac{\gamma_x^k^k}{\tau^k^k}\fprod{x_t^k-\tilde x^{k-1},x-x_{t+1}^k} \\
%&+\frac{1}{2\tau^k^k}\Big(\norm{x_t^k-x}^2-\norm{x_{t+1}^k-x}^2-\norm{x_t^k-x_{t+1}^k}^2\Big) \nonumber.
\end{align}
where the equality holds by using Lemma \ref{lem:triangle-breg}-$(d)$ and the update rule of $\hat{x}_t^k$ in line 13 of Algorithm \ref{alg:SVRPDA}.
Using  the  three-point property of  Bregman  distance  in  Lemma \ref{lem:triangle-breg}-$(b)$, the above inequality can be rewritten as 
\begin{align}\label{eq:prox-f}
f(x_{t+1}^k)-f(x)\leq & \fprod{\zeta_t,x-x_{t+1}^k}
\\&+\frac{1-\gamma_x^k}{\tau^k}\Big(\bD_X(x,x_t^k)-\bD_X(x,x_{t+1}^k)-\bD_X(x_{t+1}^k,x_t^k)\Big)\nonumber\\ \nonumber
&\eyh{+\frac{\gamma_x^k}{\tau^kT^{k-1}}\sum_{\ell=1}^{T^{k-1}}\Big(\bD_X(x,x_{\ell}^{k-1})-\bD_X(x,x_{t+1}^k)-\bD_X(x_{t+1}^k,x_{\ell}^{k-1})\Big).}
\end{align}
Recall that $\delta^x_t= \zeta_t-\grad_x\Phi(x_t^k,y_{t+1}^k)$. we add and subtract $\fprod{\grad_x\Phi(x_t^k,y_{t+1}^k),x-x_{t+1}^k}$ to the right-hand side of \eqref{eq:prox-f}, rearranging the terms, and using Definition \ref{def:mean-bregman} lead to
\begin{align}\label{eq:prox-f-2}
f(x_{t+1}^k)-f(x)\leq & \fprod{\grad_x\Phi(x_t^k,y_{t+1}^k),x-x_{t+1}^k}+\fprod{\delta^x_t,x-x_{t+1}^k} \\
& +\frac{1-\gamma_x^k}{\tau^k}\Big(\bD_X(x,x_t^k)-\bD_X(x,x_{t+1}^k)-\bD_X(x_{t+1}^k,x_t^k)\Big)\nonumber\\ \nonumber
&\eyh{+\frac{\gamma_x^k}{\tau^k}\Big(\cD_X^k(x)-\bD_X(x,x_{t+1}^k)-\cD_X^k(x_{t+1}^k)\Big).} 
\end{align}
The first inner product in the right hand side of \eqref{eq:prox-f-2} can be bounded as follows
\begin{align}\label{eq:lip-conv-x}
\nonumber&\fprod{\grad_x\Phi(x_t^k,y_{t+1}^k),x-x_{t+1}^k} \\&= \fprod{\grad_x\Phi(x_t^k,y_{t+1}^k),x-x_t^k}+\fprod{\grad_x\Phi(x_t^k,y_{t+1}^k),x_t^k-x_{t+1}^k}\nonumber \\
&\leq \Phi(x,y_{t+1}^k)-\Phi(x_t^k,y_{t+1}^k) +\fprod{\grad_x\Phi(x_t^k,y_{t+1}^k),x_t^k-x_{t+1}^k} \nonumber \\
&\leq \Phi(x,y_{t+1}^k)-\Phi(x_{t+1}^k,y_{t+1}^k) +L_{xx}\bD_X(x_{t+1}^k,x_t^k),
\end{align}
where the first inequality hold due to convexity of $\Phi(\cdot,y)$, for any $y\in\cY$; in the second inequality \eqref{eq:Lxx_bound} and the fact that $\bD_X(x,\bar{x})\geq \frac{1}{2}\norm{x-\bar{x}}^2_\cX$, for any $x,\bar{x}\in\cX$, are used. Now, we use \eqref{eq:lip-conv-x} within \eqref{eq:prox-f-2}, then rearranging the terms leads to the following result.
\begin{eqnarray}\label{eq:prox-f-acc}
\lefteqn{f(x_{t+1}^k)-f(x)-\Phi(x,y_{t+1}^k)+\Phi(x_{t+1}^k,y_{y+1}^k)\leq} \nonumber\\
&& \fprod{\delta^x_t,x-x_{t+1}^k}+L_{xx}\bD_X(x_{t+1}^k,x_t^k) \nonumber\\
&& +\frac{1-\gamma_x^k}{\tau^k}\Big(\bD_X(x,x_t^k)-\bD_X(x,x_{t+1}^k)-\bD_X(x_{t+1}^k,x_t^k)\Big)\nonumber\\
&&\eyh{+\frac{\gamma_x^k}{\tau^k}\Big(\cD_X^k(x)-\bD_X(x,x_{t+1}^k)-\cD_X^k(x_{t+1}^k)\Big).} 
\end{eqnarray}
Combining the results of \eqref{eq:prox-h-acc-3} and \eqref{eq:prox-f-acc} and using the definition of $\cL(x,y)$  \eyh{as well as Definition \ref{def:mean-bregman}} lead to the following:
\begin{align}\label{eq:lagrange-bound}
&\cL(x_{t+1}^k,y)-\cL(x,y_{t+1}^k)\leq \fprod{\delta_t^y,y_{t+1}^k-y}+\fprod{\delta^x_t,x-x_{t+1}^k}+\fprod{q_t^k-\bar{q}_t^k,y_t^k-y}\nonumber\\
&+ \fprod{\bar{q}_t^k,y_t^k-y}-\fprod{\bar{q}_{t+1}^k,y_{t+1}^k-y} +\frac{\gamma_x^k}{\tau^k}\Big(\eyh{\cD^k_X(x)}-\bD_X(x,x_{t+1}^k)\Big)\nonumber\\
&+\frac{1-\gamma_x^k}{\tau^k}\Big(\bD_X(x,x_t^k)-\bD_X(x,x_{t+1}^k)\Big) -\frac{\gamma_x^k}{\tau^k}\eyh{\cD_X^k(x_{t+1}^k)}+\frac{L_{yx}^2}{\alpha}\bD_X(x_t^k,x_{t-1}^k) \nonumber\\
&-\Big(\frac{1-\gamma_x^k}{\tau^k}-L_{xx}\Big)\bD_X(x_{t+1}^k,x_t^k)  +\frac{\gamma_y^k}{\sigma^k}\Big(\eyh{\cD_Y^k(y)} -\bD_Y(y,y_{t+1}^k)\Big)\nonumber\\
&+\frac{1-\gamma_y^k}{\sigma^k}\Big(\bD_Y(y,y_t^k)-\bD_Y(y,y_{t+1}^k)\Big)-\frac{\gamma_y^k}{\sigma^k}\eyh{\cD_Y^k(y_{t+1}^k)} \nonumber\\
&+\frac{L_{yy}^2}{\beta}\bD_Y(y_t^k,y_{t-1}^k)-\Big(\frac{1-\gamma_y^k}{\sigma^k}-(\alpha+\beta)\Big)\bD_Y(y_{t+1}^k,y_t^k) .
\end{align}
Finally, we use Lemma \ref{lem:error-inner-prod} to bound the first three inner products in the right hand side of \eqref{eq:lagrange-bound}, then we add and subtract $B_t^k$ to the right hand side which conclude the result.
\qed
\noindent Now we are ready to prove the results in Theorem \ref{thm:main} and Corollary \ref{cor:constant}, \ref{cor:nonconstant}.

{\bf Proof of Theorem \ref{thm:main}.} 
Consider the result in Lemma \ref{lem:one-step}, summing the inequality over $t=0,\hdots,T^k-1$, divide by $T^k$, and using the step-size conditions \eqref{eq:step-tau-sigma} in Assumption \ref{assum:step} one can conclude that for $k\geq 1$,
{
\begin{eqnarray}\label{eq:one-step-sum}
\lefteqn{\frac{1}{T^k}\sum_{t=0}^{T^k-1}(\cL(x_{t+1}^k,y)-\cL(x,y_{t+1}^k))\leq} \\
&&  \frac{1}{T^k}\sum_{t=0}^{T^k-1}(A_t^k-B_t^k)+\frac{1}{T^k}(\fprod{\bar{q}_0^k,y_0^k-y}-\fprod{\bar{q}_{T^k}^k,y_{T^k}^k-y})+\nonumber \\
&& \frac{1}{\eta^k T^k} \Big(\bD_X(x,u_0^k)-\bD_X(x,u_{T^k}^k)\Big)  +\frac{1}{\eta^k T^k} \Big(\bD_Y(y,v_0^k)-\bD_Y(y,v_{T^k}^k)\Big) \nonumber \\
&&+\frac{1}{\eta^k T^k} \Big(\bD_Y(y,w_0^k)-\bD_Y(y,w_{T^k}^k)\Big) +\frac{\gamma_x^k}{\tau^k}\Big(\eyh{\cD_X^k(x) - \cD_X^{k+1}(x)}\Big)\nonumber \\
&&+\frac{1-\gamma_x^k}{\tau^k T^k}\Big(\bD_X(x,x_0^k)-\bD_X(x,x_{T^k}^k)\Big) \nonumber\\
&&+\Big((6C_\cX L_{xx}^2+8C_\cY L_{yx}^2)\eta^k-\frac{\gamma_x^k}{\tau^k}\Big)\frac{1}{T^k}\sum_{t=0}^{T^k-1}\eyh{\cD_X^k(x_{t+1}^k)}\nonumber \\
&& +\frac{M_x^k}{T^k}\Big(\bD_X(x_0^k,x_{-1}^k)-\bD_X(x_{T^k}^k,x_{T^k-1}^k)\Big) \nonumber\\
%&&+\frac{\gamma_x^k}{2\tau^k}\Big(\norm{x-\tilde x^{k-1}}^2-\frac{1}{T^k}\sum_{t=0}^{T^k-1}\norm{x_{t+1}^k-x}^2\Big)+\frac{1-\gamma_x^k^k}{2\tau^k^k T^k}\Big(\norm{x_0^k-x}^2-\norm{x_{T^k}^k-x}^2\Big) \nonumber\\
%&&+\big((4L_{yx}^2+L_{xx}^2)\eta^k^k-\frac{\gamma_x^k^k}{2\tau^k^k}\big)\frac{1}{T^k}\sum_{t=0}^{T^k-1}\norm{\tilde x^{k-1}-x_{t+1}^k}^2 ~ +\frac{L_{yx}^2}{2\alpha T^k}\Big(\norm{x_0^k-x_{-1}^k}^2-\norm{x_{T^k-1}^k-x_{T^k}^k}^2\Big) \nonumber\\
&& +\frac{\gamma_y^k}{\sigma^k}\Big(\eyh{\cD_Y^k(y) -\cD_Y^{k+1}(y)}\Big)+\frac{1-\gamma_y^k}{\sigma^k T^k} \Big(\bD_Y(y,y_0^k)-\bD_Y(y,y_{T^k}^k)\Big)\nonumber\\
&&+\Big((6C_\cX L_{xy}^2+8C_\cY L_{yy}^2)\eta^k-\frac{\gamma_y^k}{\sigma^k}\Big)\frac{1}{T^k}\sum_{t=0}^{T^k-1}\eyh{\cD_Y^k(y_{t+1}^k)}\nonumber \\
&&+\frac{M_y^k}{T^k}\Big( \bD_Y(y_0^k,y_{-1}^k)-\bD_Y(y_{T^k}^k,y_{T^k-1}^k)\Big). \nonumber
%&& +\frac{\gamma_y^k^k}{2\sigma^k^k}\Big(\norm{y-\tilde y^{k-1}}^2-\frac{1}{T^k}\sum_{t=0}^{T^k-1}\norm{y_{t+1}^k-y}^2\Big)+\frac{1-\gamma_y^k^k}{2\sigma^k^k T^k}\Big(\norm{y_0^k-y}^2-\norm{y_{T^k}^k-y}^2\Big)\nonumber\\
%&&\big((4L_{xy}^2+2L_{yy}^2)\eta^k^k-\frac{\gamma_y^k^k}{2\sigma^k^k}\big)\frac{1}{T^k}\sum_{t=0}^{T-1}\norm{\tilde y^{k-1}-y_{t+1}^k}^2 ~ +\frac{L_{yy}^2}{2\beta T^k}\Big(\norm{y_0^k-y_{-1}^k}^2-\norm{y_{T^k-1}^k-y_{T^k}^k}^2\Big) .\nonumber
\end{eqnarray}}%
%Using Lemma \ref{} we can bound $\fprod{\bar{q}_{T^k}^k,y-y_{T^k}^k}$ as follows:
%\begin{align}\label{eq:q_T}
%\fprod{\bar{q}_{T^k}^k,y-y_{T^k}^k}\leq \frac{L_{yx}^2}{2\alpha}\norm{x_{T^k}^k-x_{T^k-1}^k}^2+\frac{L_{yy}^2}{2\beta}\norm{y_{T^k}^k-y_{T^k-1}^k}^2+\frac{\alpha+\beta}{2}\norm{y_{T^k}^k- y}^2 
%\end{align}

Recall that $(x_{-1}^{k+1},y_{-1}^{k+1})=(x_{T^k-1}^k,y_{T^k-1}^k)$, $(x_0^{k+1},y_0^{k+1})=(x_{T^k}^k,y_{T^k}^k)$ for any $k\geq 1$, and $(\tilde x^k,\tilde y^k)=\frac{1}{T^k}\sum_{t=0}^{T^k-1}(x_{t+1}^k,y_{t+1}^k)$, then using \eqref{eq:step-gamma} in Assumption \ref{assum:step} and Jensen's inequality, i.e., $\frac{1}{T}\sum_{t=0}^{T-1}\phi(z_t)\geq \phi\left(\tfrac{1}{T}\sum_{t=0}^{T-1}z_t\right)$ for any $T>0$ and any convex function $\phi$, for the convex functions $\cL(\cdot,y)$ and $-\cL(x,\cdot)$, 
%$\bD_X(x,\cdot)$, and $\bD_Y(y,\cdot)$ (see Assumptions \ref{assum:lip} and \ref{assum:bregman}-$(b)$) 
we obtain the following %result
{
\begin{eqnarray}\label{eq:lagrange-average}
\lefteqn{\cL(\tilde x^k,y)-\cL(x,\tilde y^k)\leq} \\
&&  \frac{1}{T^k}\sum_{t=0}^{T^k-1}(A_t^k-B_t^k)+\frac{1}{T^k}(\fprod{\bar{q}_0^k,y_0^k-y}-\fprod{\bar{q}_{0}^{k+1},y_{0}^{k+1}-y}) 
\nonumber \\
&&+\frac{1}{\eta^k T^k} \Big(\bD_X(x,u_0^k)-\bD_X(x,u_{0}^{k+1})\Big) \nonumber \\
&& +\frac{1}{\eta^k T^k} \Big(\bD_Y(y,v_0^k)-\bD_Y(y,v_{0}^{k+1})\Big) +\frac{1}{\eta^k T^k} \Big(\bD_Y(y,w_0^k)-\bD_Y(y,w_{0}^{k+1})\Big) \nonumber \\
&&+\frac{\gamma_x^k}{\tau^k}\Big(\eyh{\cD_X^k(x) - \cD_X^{k+1}(x)}\Big)+\frac{1-\gamma_x^k}{\tau^k T^k}\Big(\bD_X(x,x_0^k)-\bD_X(x,x_{0}^{k+1})\Big) \nonumber\\
&&+\frac{M_x^k}{T^k}\Big(\bD_X(x_0^k,x_{-1}^k)-\bD_X(x_{0}^{k+1},x_{-1}^{k+1})\Big) \nonumber\\
&& +\frac{\gamma_y^k}{\sigma^k}\Big(\eyh{\cD_Y^k(y) -\cD_Y^{k+1}(y)}\Big)+\frac{1-\gamma_y^k}{\sigma^k T^k} \Big(\bD_Y(y,y_0^k)-\bD_Y(y,y_{0}^{k+1})\Big)\nonumber\\
&&+\frac{M_y^k}{T^k}\Big( \bD_Y(y_0^k,y_{-1}^k)-\bD_Y(y_{0}^{k+1},y_{-1}^{k+1})\Big). \nonumber
%+\frac{1}{2\eta^kT^k}\Big(\norm{x-u_0^k}^2-\norm{x-u_{0}^{k+1}}^2\Big) \nonumber \\
% && +\frac{1}{2\eta^kT^k}\Big(\norm{y-v_0^k}^2-\norm{y-v_{0}^{k+1}}^2\Big) +\frac{1}{2\eta^kT^k}\Big(\norm{y-w_0^k}^2-\norm{y-w_{0}^{k+1}}^2\Big) \nonumber\\
% &&+\frac{\gamma_x^k}{2\tau^k}\Big(\norm{x-\tilde x^{k-1}}^2-\norm{\tilde{x}^{k}-x}^2\Big)+\frac{1-\gamma_x^k}{2\tau^k T^k}\Big(\norm{x_0^k-x}^2-\norm{x_{0}^{k+1}-x}^2\Big) \nonumber\\
% && +\frac{L_{yx}^2}{2\alpha T^k}\Big(\norm{x_0^k-x_{-1}^k}^2-\norm{x_{-1}^{k+1}-x_{0}^{k+1}}^2\Big) \nonumber\\
% && +\frac{\gamma_y^k}{2\sigma^k}\Big(\norm{y-\tilde y^{k-1}}^2-\norm{\tilde{y}^{k+1}-y}^2\Big)+\frac{1-\gamma_y^k}{2\sigma^k T^k}\Big(\norm{y_0^k-y}^2-\norm{y_{0}^{k+1}-y}^2\Big)\nonumber\\
% && +\frac{L_{yy}^2}{2\beta T^k}\Big(\norm{y_0^k-y_{-1}^k}^2-\norm{y_{-1}^{k+1}-y_{0}^{k+1}}^2\Big) .\nonumber
\end{eqnarray}}%

Now multiplying \eqref{eq:lagrange-average} by $T^{k}$, leads to %using the fact that $\frac{\gamma^k_x}{\tau^k}T^k=c_1(k+1), \frac{1-\gamma^k_x}{\tau^k}=c_2(k+1), \frac{\gamma^k_y}{\sigma^k}T^k=c_3(k+1), \frac{1-\gamma^k_y}{\sigma^k}=c_4(k+1)$, for some $c_1,c_2,c_3,c_4\in\reals_{++}$, and $\eta^k=(\eta (k+1))^{-1}$, %$\frac{s^{k-1}}{s^k}=\frac{T^{k-1}}{T^k}\geq \max\big\{\frac{\eta^{k-1} T^{k-1}}{\eta^k T^k}, \frac{\gamma_x^k\tau^{k-1}}{\gamma_x^{k-1}\tau^k},\frac{(1-\gamma_x^{k})\tau^{k-1}T^{k-1}}{(1-\gamma_x^{k-1})\tau^{k}T^{k}},\frac{\gamma_y^k\sigma^{k-1}}{\gamma_y^{k-1}\sigma^k},$ $\frac{(1-\gamma_y^{k})\sigma^{k-1}T^{k-1}}{(1-\gamma_y^{k-1})\sigma^{k}T^{k}} \big\}$, 
%we obtain the following:
{
\begin{align}\label{eq:lagrange-rho}
&T^k(\cL(\tilde x^k,y)-\cL(x,\tilde y^k))\\
&\leq \sum_{t=0}^{T^k-1}(A_t^k-B_t^k)+\fprod{\bar{q}_0^k,y_0^k-y}-\fprod{\bar{q}_0^{k+1},y_0^{k+1}-y}+\frac{1}{\eta^k}(Q^k-Q^{k+1})\nonumber\\
&+\frac{\gamma_x^k}{\tau^k}T^k\Big(\eyh{\cD_X^k(x) - \cD_X^{k+1}(x)}\Big)+\frac{1-\gamma_x^k}{\tau^k}\Big(\bD_X(x,x_0^k)-\bD_X(x,x_{0}^{k+1})\Big) \nonumber\\
&+M_x^k\Big(\bD_X(x_0^k,x_{-1}^k)-\bD_X(x_{0}^{k+1},x_{-1}^{k+1})\Big) \nonumber\\
& +\frac{\gamma_y^k}{\sigma^k}T^k\Big(\eyh{\cD_Y^k(y) - \cD_Y^{k+1}(y)}\Big)+\frac{1-\gamma_y^k}{\sigma^k} \Big(\bD_Y(y,y_0^k)-\bD_Y(y,y_{0}^{k+1})\Big)\nonumber\\
&+M_y^k\Big( \bD_Y(y_0^k,y_{-1}^k)-\bD_Y(y_{0}^{k+1},y_{-1}^{k+1})\Big),\nonumber
\end{align}}
% \begin{eqnarray}\label{eq:lagrange-rho}
% \lefteqn{T^k(\cL(\tilde x^k,y)-\cL(x,\tilde y^k))\leq \sum_{t=0}^{T^k-1}(A_t^k-B_t^k)} \\
% &&+\fprod{\bar{q}_0^k,y_0^k-y}+ \frac{1}{2\eta^{k}}\Big(\norm{x-u_0^k}^2+\norm{y-v_0^k}^2+\norm{y-w_0^k}^2\Big) \nonumber\\
% &&+\frac{\gamma_x^{k}}{2\tau^{k}}T^k\norm{x-\tilde x^{k-1}}^2+\frac{1-\gamma_x^{k}}{2\tau^{k}}\norm{x_0^k-x}^2 +\frac{L_{yx}^2}{2\alpha}\norm{x_0^k-x_{-1}^k}^2 \nonumber\\
% && +\frac{\gamma_y^{k}}{2\sigma^{k}}T^k\norm{y-\tilde y^{k-1}}^2+\frac{1-\gamma_y^{k}}{2\sigma^{k}}\norm{y_0^k-y}^2 +\frac{L_{yy}^2}{2\beta} \norm{y_0^k-y_{-1}^k}^2 \nonumber\\
% &&-\Big[\fprod{\bar{q}_0^{k+1},y_0^{k+1}-y}+ \frac{1}{2\eta^{k}}\Big(\norm{x-u_0^{k+1}}^2+\norm{y-v_0^{k+1}}^2+\norm{y-w_0^{k+1}}^2\Big) \nonumber\\
% &&+\frac{\gamma_x^{k}}{2\tau^{k}}T^k\norm{x-\tilde x^{k}}^2+\frac{1-\gamma_x^{k}}{2\tau^{k}}\norm{x_0^{k+1}-x}^2 +\frac{L_{yx}^2}{2\alpha} \norm{x_0^{k+1}-x_{-1}^{k+1}}^2 \nonumber\\
% && +\frac{\gamma_y^{k}}{2\sigma^{k}}T^k\norm{y-\tilde y^{k}}^2+\frac{1-\gamma_y^{k}}{2\sigma^{k}}\norm{y_0^{k+1}-y}^2 +\frac{L_{yy}^2}{2\beta} \norm{y_0^{k+1}-y_{-1}^{k+1}}^2\Big] . \nonumber
% \end{eqnarray}
where %$C_1^k\triangleq \diag([\frac{\gamma_x^k}{\tau^k}T^k\id_n,~\frac{\gamma_y^k}{\sigma^k}T^k\id_m])$, $C_2^k\triangleq \diag([\frac{1-\gamma_x^k}{\tau^k}\id_n,\frac{1-\gamma_y^k}{\sigma^k}\id_m])$, and
$Q^k\triangleq  \bD_X(x,u_0^k)+\bD_Y(x,v_0^k)+\bD_Y(x,w_0^k)$ for any $k\geq 1$. 
%Now we consider two cases whether the Bregman diameters $B_X$ and $B_Y$ are bounded. 
Next, we consider two cases depending on selecting the step-sizes.

{\bf Part I)} In this scenario, we consider constant step-sizes and parameters, i.e., $T^k=\bar{T}$, $\tau^k=\tau$, $\sigma^k=\sigma$, $\gamma_x^k=\gamma_x$, $\gamma_y^k=\gamma_y$, $\eta^k=\eta$. Then, summing \eqref{eq:lagrange-rho} over $k$ from 1 to $K$, using \eqref{eq:bound-telescop} and Jensen's inequality we conclude  that %for $K\geq 1$, 
{
\begin{eqnarray}\label{eq:lagrange-rho-sum-constant}
\lefteqn{\bar{T}K(\cL(\tilde x^{(K)},y)-\cL(x,\tilde y^{(K)}))\leq }\nonumber\\
&&\sum_{k=1}^K\sum_{t=0}^{\bar{T}-1}(A_t^k-B_t^k)~-\fprod{\bar{q}_0^{K+1},y^{K+1}_0-y}+\frac{1}{\eta}(Q^0-Q^{K+1}) \nonumber\\
&&+\frac{\gamma_x}{\tau}\bar{T}(\eyh{\cD_X^1(x) - \cD_X^{K+1}(x)})+\frac{1-\gamma_x}{\tau}(\bD_X(x,x_0^{1})-\bD_X(x,x_0^{K+1}))\nonumber\\
&&+M_x(\bD_X(x_0^{1},x_{-1}^{1})-\bD_X(x_0^{K+1},x_{-1}^{K+1}))\nonumber\\
&& +\frac{\gamma_y}{\sigma}\bar{T}(\eyh{\cD_Y^1(y) - \cD_Y^{K+1}(y)})+\frac{1-\gamma_y}{\sigma}(\bD_Y(y,y_0^{1})-\bD_Y(y,y_0^{K+1}))\nonumber\\
&&+M_y(\bD_Y(y_0^{1},y_{-1}^{1})-\bD_Y(y_0^{K+1},y_{-1}^{K+1})).
\end{eqnarray}}%
Note that the inner product term $\fprod{\bar{q}_0^{K+1},y^{K+1}_0-y}$ in the right hand side of \eqref{eq:lagrange-rho-sum} can be lower bounded using \eqref{eq:inner-prod-sum} as follows
{
\begin{align}\label{eq:lower-inner}
\fprod{\bar{q}_0^{K+1},y^{K+1}_0-y}\nonumber\geq & -\frac{L_{yx}^2}{\alpha}\bD_X(x_0^{K+1},x_{-1}^{K+1})-\frac{L_{yy}^2}{\beta}\bD_Y(y_0^{K+1},y_{-1}^{K+1})\\ &-(\alpha+\beta)\bD_Y(y,y_0^{K+1}).
\end{align}}%
Using \eqref{eq:lower-inner} within \eqref{eq:lagrange-rho-sum-constant}, the step-size condition \eqref{eq:step-tau-sigma} in Assumption \ref{assum:step}, and then dropping the nonpositive terms lead to
{
\begin{eqnarray}\label{eq:final-lagrange-constant}
\lefteqn{\bar{T}K(\cL(\tilde x^{(K)},y)-\cL(x,\tilde y^{(K)}))\leq }\\
&&\sum_{k=1}^K\sum_{t=0}^{\bar{T}-1}(A_t^k-B_t^k)+\frac{1}{\eta}Q^0+\frac{\gamma_x}{\tau}\bar{T}\eyh{\cD_X^1(x)}+\frac{1-\gamma_x}{\tau}\bD_X(x,x_0^{1})\nonumber\\
&&+M_x\bD_X(x_0^{1},x_{-1}^{1})+\frac{\gamma_y}{\sigma}\bar{T}\eyh{\cD_Y^1(y)}+\frac{1-\gamma_y}{\sigma}\bD_Y(y,y_0^{1})+M_y\bD_Y(y_0^{1},y_{-1}^{1}).\nonumber
\end{eqnarray}}%
Finally, taking supremum over $z=(x,y)\in Z=X\times Y$, then taking expectation $\mE[\cdot]$ and using the fact that $\mE[A_t^k]\leq B_t^k$ from Lemma \ref{lem:upper-noise} leads to the result in \eqref{eq:thm-main-result-constant}.

{\bf Part II)}
Let $\{s^k\}_{k\geq 1}\subset\reals_+$ be a sequence and $\{P^k\}_{k\geq 1}\subset\reals_+$ be a bounded sequence, i.e., there exists $\Delta>0$ such that $\sup_{k\geq 1}P^k\leq \Delta$, then %one can observe that
{
\begin{align}\label{eq:bound-telescop}
 \sum_{k=1}^K s^k(P^k-P^{k+1})=s^1P^1+\Big[\sum_{k=2}^K(s^k-s^{k-1})\Delta\Big]-s^KP^{K+1}\leq s^K(\Delta-P^{K+1}) . 
\end{align}}%
Therefore, assuming that the Bregman diameter is bounded, %and the fact that $\{\frac{1}{\eta^k}\}_{k\geq 1}$, $\{\frac{\gamma_x^k}{\tau^k}T^k,\frac{1-\gamma_x^k}{\tau^k},M_x^k\}_{k\geq 1}$,  $\{\frac{\gamma_y^k}{\sigma^k}T^k,\frac{1-\gamma_y^k}{\sigma^k},M_y^k\}_{k\geq 1}\subset \reals_+$ are nondecreasing sequences-- see Assumption \ref{assum:step} and Remark \ref{rem:step} for a particular choice of parameters, 
one can use \eqref{eq:bound-telescop} for each term in the right hand side of \eqref{eq:lagrange-rho} involving differences of two consecutive Bregman distance functions. % when summing over $k$ from 1 to $K$, leading to  
% \begin{align}\label{eq:sum-R-k}
% \sum_{k=1}^K T^k(R^{k-1}-R^{k}) & = T^1R^0+\Big[\sum_{k=2}^{K} (T^k-T^{k-1})R^{k-1}\Big]-T^KR^{K}\nonumber\\
% &\leq T^1R^0+(T^K-T^1)\Big(\frac{\gamma_x}{\tau}B_x+\frac{\gamma_y}{\sigma}B_y\Big)\leq T^K\Big(\frac{\gamma_x}{\tau}B_x+\frac{\gamma_y}{\sigma}B_y\Big).
% \end{align} 
Hence, summing \eqref{eq:lagrange-rho} over $k$ from 1 to $K$, using \eqref{eq:bound-telescop} and Jensen's inequality we conclude  that for $K\geq 1$
{
\begin{eqnarray}\label{eq:lagrange-rho-sum}
\lefteqn{S^K(\cL(\tilde x^{(K)},y)-\cL(x,\tilde y^{(K)}))\leq }\nonumber\\
&&\sum_{k=1}^K\sum_{t=0}^{T^k-1}(A_t^k-B_t^k)~-\fprod{\bar{q}_0^{K+1},y^{K+1}_0-y}+\frac{1}{\eta^K}(B_X+2B_Y-Q^{K+1}) \nonumber\\
&&+\frac{\gamma_x^K}{\tau^K}T^K(B_X-\eyh{\cD_X^{K+1}(x)})+\frac{1-\gamma_x^K}{\tau^K}(B_X-\bD_X(x,x_0^{K+1}))\nonumber\\
&&+M_x^K(B_X-\bD_X(x_0^{K+1},x_{-1}^{K+1}))\nonumber\\
&& +\frac{\gamma_y^K}{\sigma^K}T^K(B_Y-\eyh{\cD_Y^{K+1}(y)})+\frac{1-\gamma_y^K}{\sigma^K}(B_Y-\bD_Y(y,y_0^{K+1}))\nonumber\\
&&+M_y^K(B_Y-\bD_Y(y_0^{K+1},y_{-1}^{K+1})),
\end{eqnarray}}%
where $S^K=\sum_{k=1}^K T^k$. Similar to the previous case the inner product can be bounded by \eqref{eq:lower-inner} and using the step-size condition \eqref{eq:step-tau-sigma} we obtain,
{
\begin{eqnarray}\label{eq:lagrange-drop-bound}
\lefteqn{S^K(\cL(\tilde x^{(K)},y)-\cL(x,\tilde y^{(K)}))\leq }\\
&&\sum_{k=1}^K\sum_{t=0}^{T^k-1}(A_t^k-B_t^k)~+\Big(\frac{1}{\eta^K}+\frac{\gamma_x^K}{\tau^K}T^K+\frac{1-\gamma_x^K}{\tau^K}+M_x^K\Big)B_X\nonumber\\
&&+\Big(\frac{2}{\eta^K}+\frac{\gamma_y^K}{\sigma^K}T^K+\frac{1-\gamma_y^K}{\sigma^K}+M_y^K\Big)B_Y .\nonumber
\end{eqnarray}}%
%where the last inequality holds due to the step-size condition \eqref{eq:step-tau-sigma} in Assumption \ref{assum:step}. Therefore, one can drop the $Q^{K+1}$ in \eqref{eq:lagrange-rho-sum}, 
Finally, with a similar argument as in Part I and using Lemma \ref{lem:upper-noise} the result in \eqref{eq:thm-main-result} can be concluded. \qed
%Next, one can take supremum over $z= (x,y)\in Z=X\times Y$, and then take expectation $\mE[\cdot]$ to obtain the following result:
% \begin{eqnarray}\label{eq:final-bound}
% \lefteqn{\mE\big[\sup_{z\in\cZ}\big(S^K(\cL(\tilde x^{(K)},y)-\cL(x,\tilde y^{(K)}))\big)\big]}\\
% &&\leq \sum_{k=1}^K\sum_{t=0}^{T^k-1}\mE[A_t^k-B_t^k]\nonumber\\
% &&+\Big(\frac{1}{\eta^K}+\frac{\gamma_x^K}{\tau^K}T^K+\frac{1-\gamma_x^K}{\tau^K}+M_x^K\Big)B_X+\Big(\frac{2}{\eta^K}+\frac{\gamma_y^K}{\sigma^K}T^K+\frac{1-\gamma_y^K}{\sigma^K}+M_y^K\Big)B_Y \nonumber\\ \nonumber
% &&\leq \Big(\frac{1}{\eta^K}+\frac{\gamma_x^K}{\tau^K}T^K+\frac{1-\gamma_x^K}{\tau^K}+M_x^K\Big)B_X+\Big(\frac{2}{\eta^K}+\frac{\gamma_y^K}{\sigma^K}T^K+\frac{1-\gamma_y^K}{\sigma^K}+M_y^K\Big)B_Y
% \end{eqnarray} 
%where in the last inequality we used $\mE[A_t^k]\leq B_t^k$ which is implied by Lemma \ref{lem:upper-noise}. 
%Dividing both sides of \eqref{eq:final-bound} by $S^K$ gives the result in \eqref{eq:thm-main-result}. \qed

{\bf Proof of Corollary \ref{cor:constant}}
Assume that the step-sizes and parameters are selected as in \eqref{eq:specific-step-const}. It is easy to verify that $\max\{\frac{1}{\eta}, \frac{\gamma_x}{\tau}\bar{T}, \frac{1-\gamma_y}{\tau}, \frac{\gamma_y}{\sigma}\bar{T},\frac{1-\gamma_y}{\sigma},M_x, M_y\}=\cO(\sqrt{n})$ and $\bar{T}=\cO(n)$; hence, we conclude that the right side of \eqref{eq:thm-main-result-constant} has the rate of $\cO(\frac{1}{K\sqrt{n}})$. Hence, the total number of gradients to achieve $\epsilon$-gap is $nK+\sum_{k=1}^KT^k=\eyh{\cO(nK)}=\cO(\frac{\sqrt{n}}{\epsilon})$. \qed

{\bf Proof of Corollary \ref{cor:nonconstant}}
Assume that the step-sizes and parameters are selected as in \eqref{eq:specific-step}, then from the fact that $T^k=T(k+1)^2$ we have that %for $ K\geq 1$,
$
S^K=\sum_{k=1}^K T(k+1)^2 = T\Big(\frac{(K+1)^3}{3}+\frac{(K+1)^2}{2}+\frac{(K+1)}{6}-1\Big)\geq \frac{T(K+1)^3}{3}$.
Moreover, $\tau^k=\cO(\frac{1}{k})$, $\sigma^k=\cO(\frac{1}{k})$, $\gamma_x^k=\cO(\frac{1}{k^2})$, $\gamma_y^k=\cO(\frac{1}{k^2})$, and $\eta^k=\cO(\frac{1}{k})$ which implies that
$\max\{\frac{1}{\eta^K}, \frac{\gamma_x^K}{\tau^K}T^K, \frac{1-\gamma_y^K}{\tau^K}, \frac{\gamma_y^K}{\sigma^K}T^K,\frac{1-\gamma_y^K}{\sigma^K},M_x^K, M_y^K\}=\cO(K)$; hence, 
we conclude that the right side of \eqref{eq:thm-main-result} has the rate of $\cO(1/K^2)$. Hence, the total number of gradients is $nK+\sum_{k=1}^KT^k=nK+S^K=\cO(\frac{n}{\sqrt{\epsilon}}+\frac{1}{\epsilon^{1.5}})$.\qed

\section{Numerical Experiment}\label{sec:num}
In this section, we implement SVR-APD with constant (SVR-APD-I) and non-constant (SVR-APD-II) step-sizes for solving DRO problem \eyh{\eqref{eq:DRO2}} described in Section \ref{sec:application}, and compare them with the state-of-the-art first-order methods designed for solving large-scale convex-concave SP problems, e.g., Stochastic Mirror Descent (SMD) \cite{juditsky2011solving} and Stochastic Mirror-prox (SMP) \cite{nemirovski2009robust}. 

%\subsection{Distributed robust optimization.}
Similar to the setup in \cite{namkoong2016stochastic}, we consider $\{a_i\}_{i=1}^n\subset \reals^m$ to be a set of features with labels $\{b_i\}_{i=1}^n\subset\{-1,+1\}^n$. We consider the logistic loss function, i.e., $\ell_i(x)=\log(1+\exp(-b_ia_i^\top x))$, and Chi-square divergence measure $V(y,\frac{1}{n}\ones_n)=\frac{1}{2}\norm{ny-\ones_n}^2$, and we set $X=[-10,10]^d$ and $\rho=50$. \eyh{Our goal is to compare the performance of the methods when $n$ (number of samples) is large. Different datasets have been used for the experiment and summary of the information can be found in Table \ref{tab:1}.} 

\begin{figure}[htb]
    \vspace*{-5mm}
    \centering
    \includegraphics[scale=0.17]{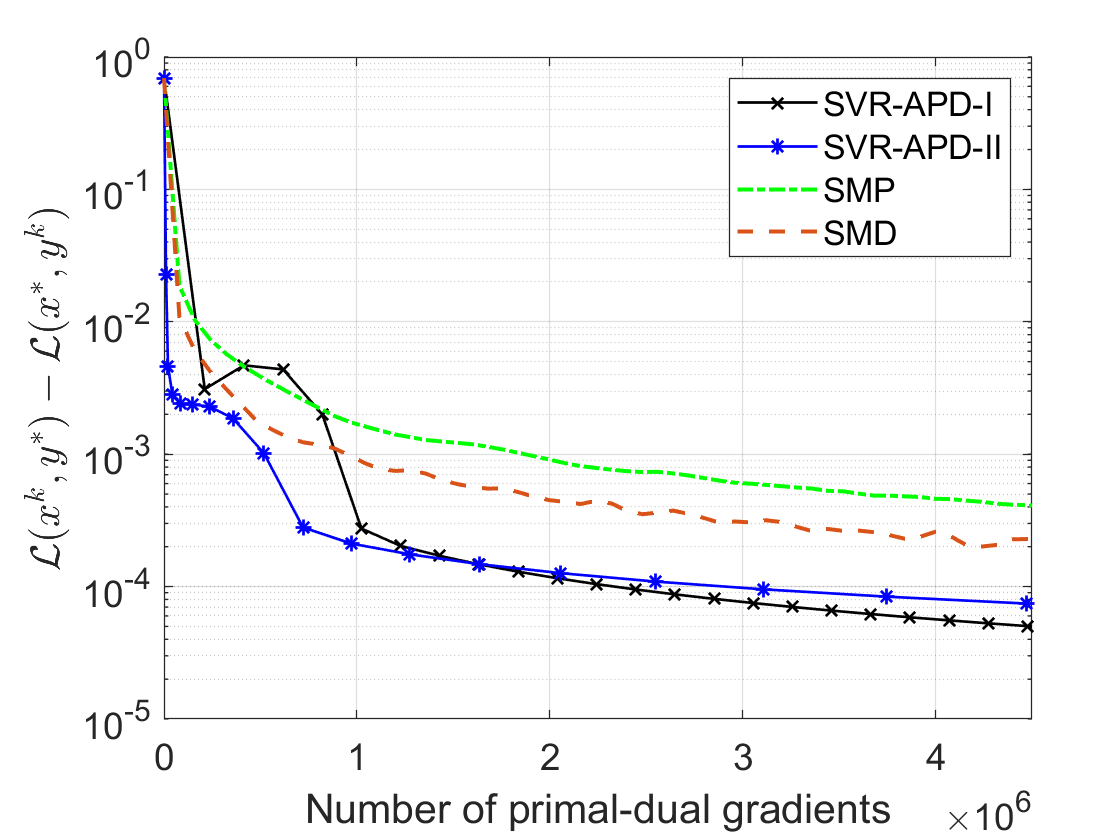} \hspace{-3mm}
        \includegraphics[scale=0.17]{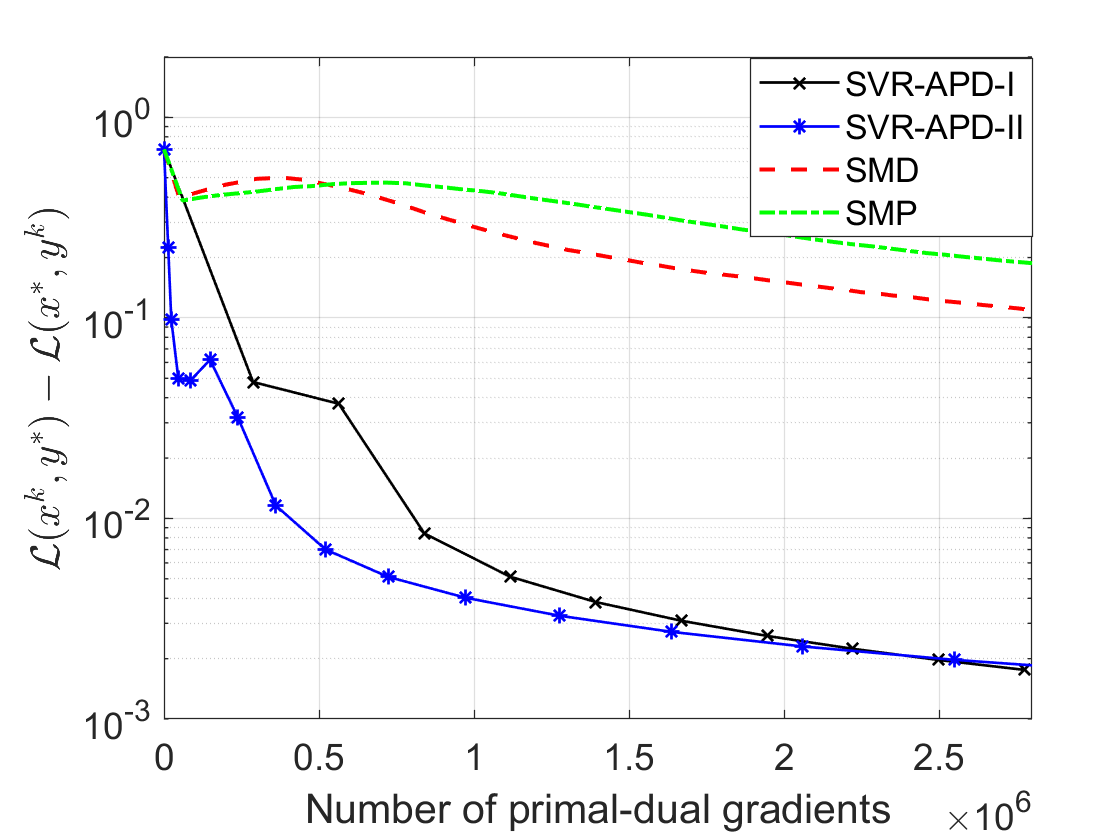}   \hspace{-3mm} \includegraphics[scale=0.17]{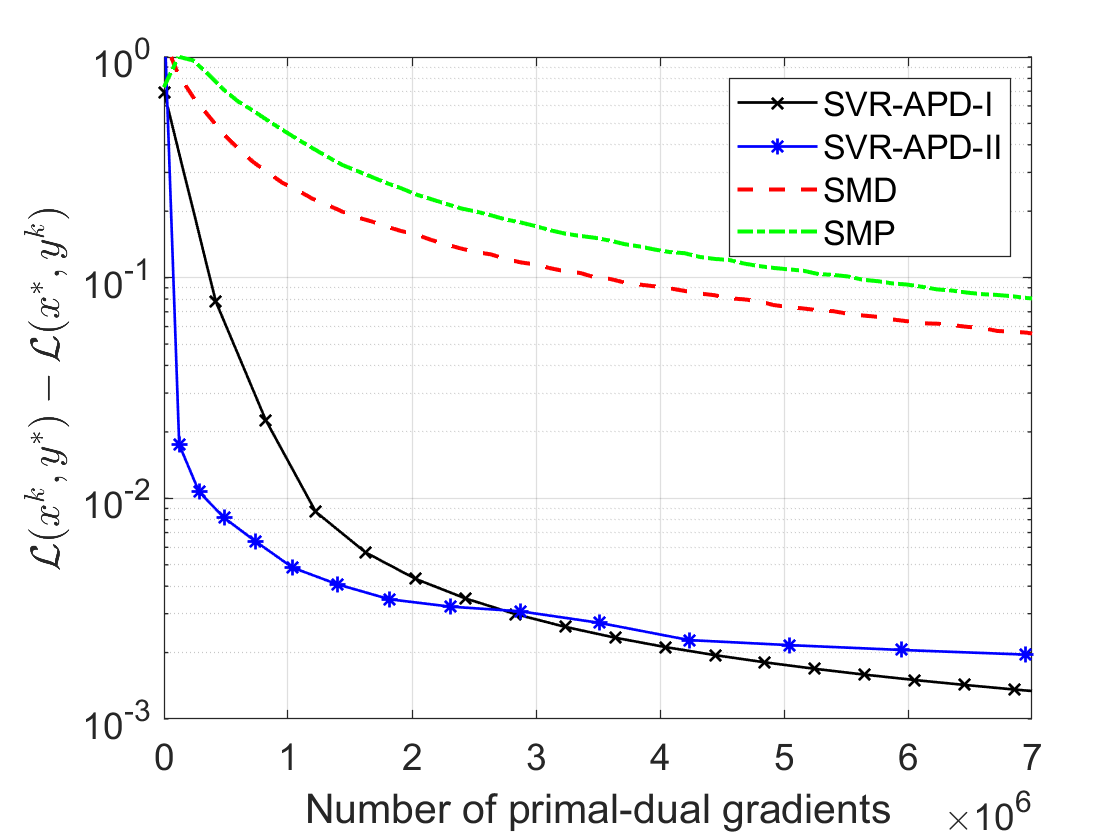}
         \includegraphics[scale=0.17]{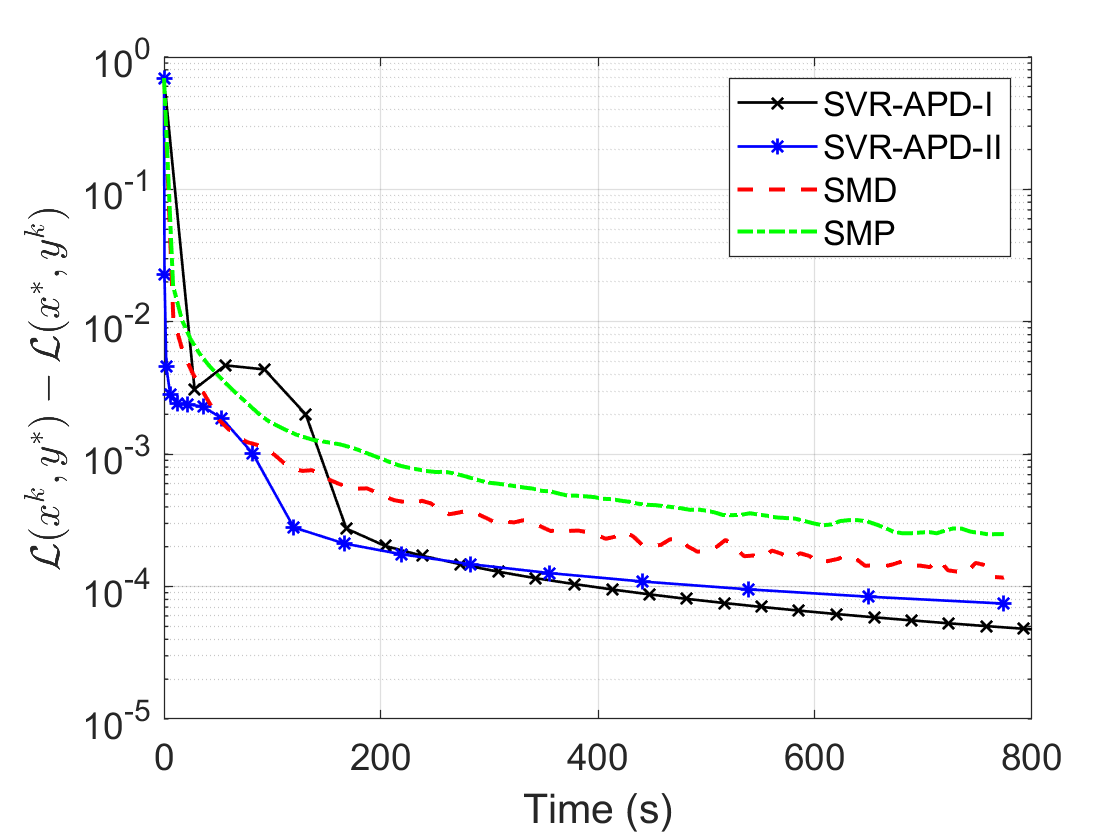} \hspace{-3mm}
        \includegraphics[scale=0.17]{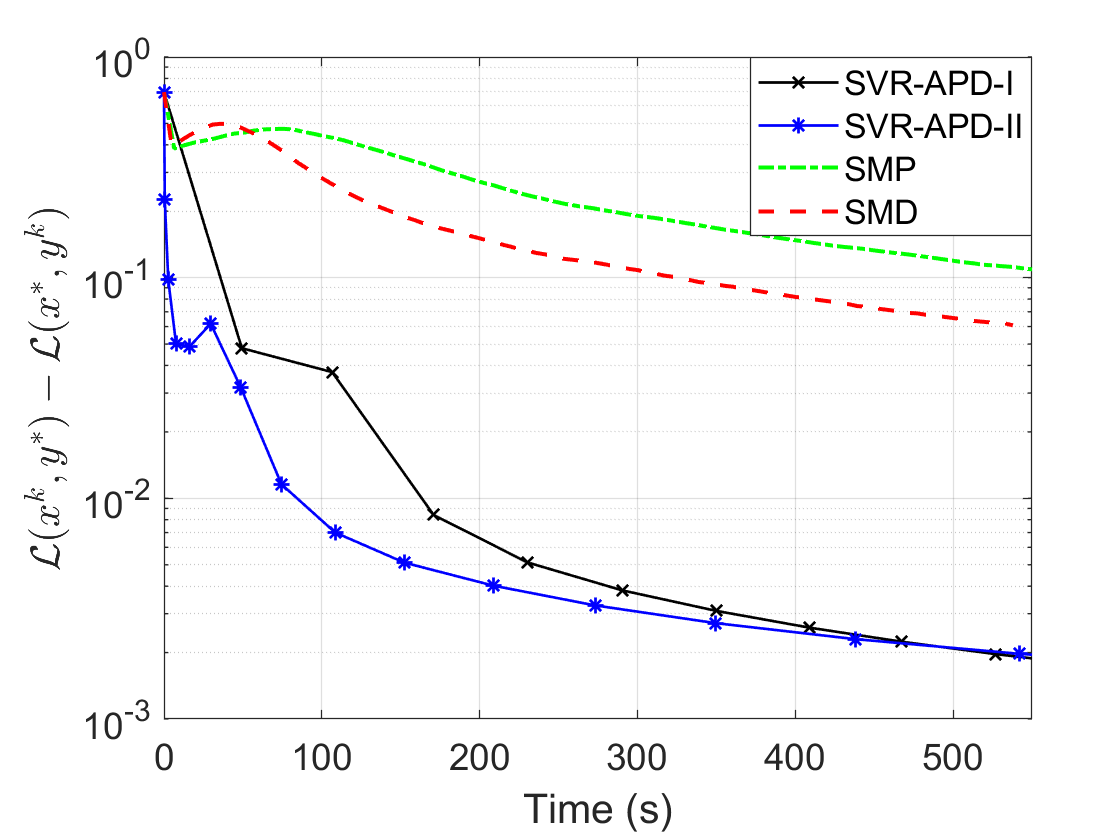}   \hspace{-3mm} \includegraphics[scale=0.17]{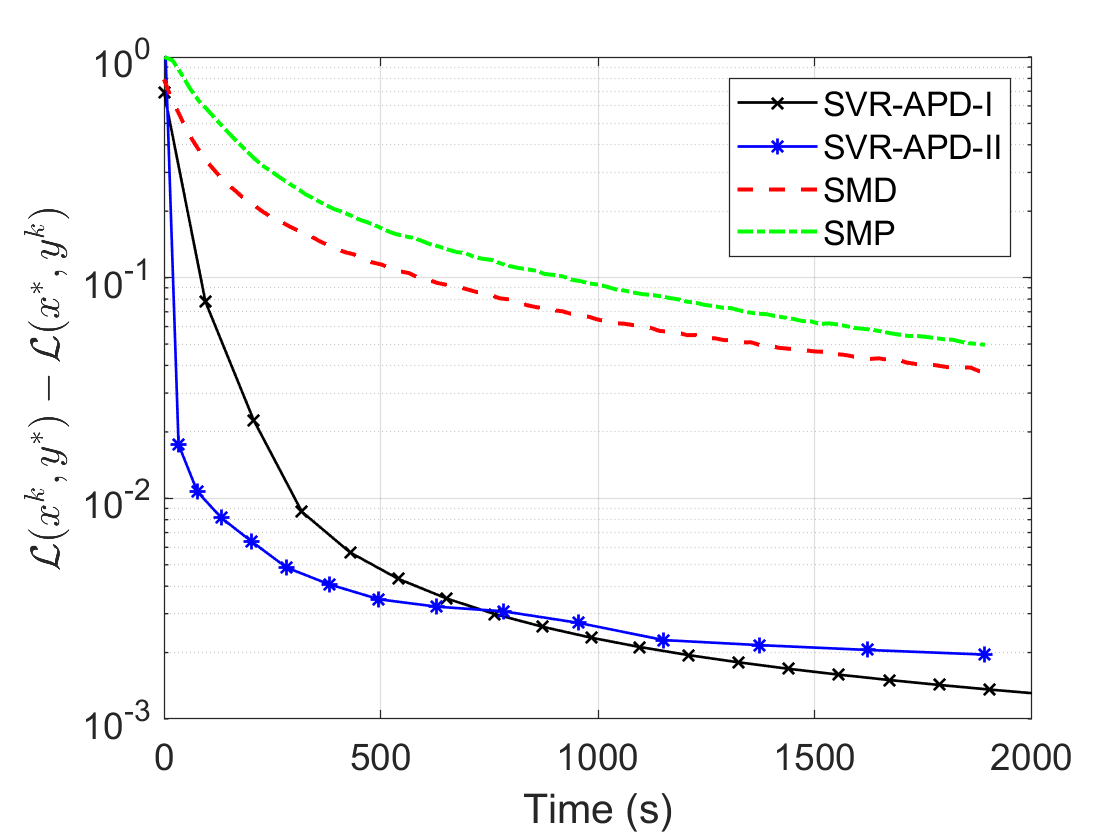}
    \caption{Comparison of the methods in terms of number of primal-dual gradients (Top row) and  time (Bottom row) for datasets Mushrooms, Phishing and a7a from left to right.}
    \label{fig:1}
\end{figure}

%We implement our proposed algorithm with constant step-size (SVR-APD-I) and non-constant step-size (SVR-APD-II) for solving \eqref{eq:DRO2} and compare with SMD and SMP. 
\eyh{To compute the projection onto the simplex-set constraint, $y\in \Delta_n$, in problem \eqref{eq:DRO2}, we choose 
%Considering problem \eqref{eq:DRO2}, there is a simplex constraint, $y\in\Delta_n$, in the maximization which has a closed form projection by considering 
the entropy Bregman distance generating function $\psi_\cY(y)=\sum_{i=1}^n y_ilog(y_i)$, where $y=[y_i]_{i=1}^n\in\reals^n$--see \cite{nesterov2005smooth} for more details. We also choose the step-sizes $\tau^k$ and $\sigma^k$ as in \eqref{eq:specific-step-const-ts} and \eqref{eq:specific-step-ts} for SVR-APD-I and SVR-APD-II, respectively and we select the parameters $\bar \gamma_x=\bar \gamma_y=1/\max\{L_x,L_y\}$.}
%Moreover, we select $T^k=10n$ and $T^k=10^3k^2$ for constant and non-constant step-sizes.
%We choose the step-size  from $10^{-i}$, for $i\in\{0,\hdots,5\}$, based on the best performance of the algorithms. For our method with constant paramters, we set $T^k=10n$, $\tau=\frac{10^{-i}}{\sqrt{n}}$, $\sigma=\frac{\tau}{10^2}$, $\gamma_x=\gamma_y=\frac{\tau}{10^2\sqrt{n}}$, and for the non-constant case $T^k=10^3 k^2$, $\tau^k=\frac{10^{-i}}{k}$, $\sigma^k=\frac{\tau^k}{10^2}$, $\gamma_x^k=\gamma_y^k=\frac{\tau^k}{10^2k}$. 
We plot the results in terms of \eyh{the difference of Lagrangian functions, i.e., $\cL(x^k,y^*)-\cL(x^*,y^k)$,} versus the number of primal-dual gradients and running time of the algorithms in Figure \ref{fig:1}. For all three experiments our method outperforms the other two schemes \eyh{and the superiority is more evident as the number of samples are getting larger}. Moreover, consistent with our results, SVR-APD-I with $\cO(\frac{\sqrt{n}}{\epsilon})$ oracle complexity  eventually outperform SVR-APD-II with $\cO(\frac{1}{\epsilon^{1.5}})$ for higher accuracy.

\begin{table}[htb]
    \centering
    \begin{tabular}{|c|c|c|c|}\hline
    &Mushrooms&Phishing&a7a\\ \hline
       \# of samples  &8124 &11055 &16100 \\ \hline
        \# of features &112 &64 &122\\ \hline
    \end{tabular}
    \caption{Datasets from UCI: http://archive.ics.uci.edu/ml/index.php}
    \label{tab:1}
\end{table}

\section{Appendix}\label{append}
The following lemma provides some fundamental properties of associate with the Bregman distance functions--see \cite{Tseng08_1J,banerjee2005clustering} for the proofs. 
\begin{lemma}\label{lem:triangle-breg}
Let $(\cU,\norm{\cdot}_\cU)$ be a finite-dimensional normed vector space with the dual space $(\cU^*,\norm{\cdot}_{\cU^*})$, $f:\cU \to \reals\cup\{+\infty\}$ be a closed convex function, $U\subset \cU$ is a closed convex set, $\psi:\cU\to\reals$ be a distance generating function which is continuously differentiable on an open set containing $\dom f$ and 1-strongly convex with respect to $\norm{\cdot}_\cU$, and $\bD_U:U\times(U\cap \dom f)\to \reals$ be a Bregman distance function associated with $\psi$. Then, the following result holds.
\begin{enumerate}[label=\alph*)]
\item Given $\bar{x}\in U\cap \dom f$, $s\in\cU$ and $t>0$, let
$x^+=\argmin_{x\in U} f(x)-\fprod{s,x}+t \bD_U(x,\bar{x})$. 
Then for all $x\in U$, the following inequality holds:
\begin{eqnarray}
f(x)\geq f(x^+)+\fprod{s,x-x^+}+ t\fprod{\grad\psi(\bar{x})-\grad\psi(x^+),x-x^+}. \label{eq_app:bregman}
\end{eqnarray}

\item For all $x\in U$ and $y,z\in U\cap\dom f$, 
$
    \fprod{\grad\psi(z)-\grad\psi(y),x-z} = \bD_U(x,y)-\bD_U(x,z)-\bD_U(z,y).
$

\item Given the update of $x^+$ in (a), for all $x\in U$ the following inequality holds:
\begin{eqnarray}
 f(x^+)-f(x)+\fprod{s,x-\bar{x}}\leq t\big(\bD_U(x,\bar{x})-\bD_U(x,x^+)\big)+\frac{1}{2t}\norm{s}_{\cU^*}^2. \label{eq_app:bregman-aux}
\end{eqnarray}

\item Assuming $\psi(\cdot)$ is a closed function, then 
$\grad\psi^*(\grad\psi(x)) = x$, for all  $x\in \dom\grad\psi\subset\cU,$ and $\grad\psi(\grad\psi^*(y)) = y$, for all $y\in \dom\grad\psi^*\subset\cU^*$.
\end{enumerate}
\end{lemma}
\bibliographystyle{spmpsci_unsrt} 
\bibliography{papers}

\end{document}